\documentclass[11pt]{amsart}
\usepackage[letterpaper,margin=1in]{geometry}                
\usepackage{amssymb,amsmath,amsfonts}
\usepackage{xcolor}
\usepackage[colorlinks,
             linkcolor=blue,
             citecolor=green!70!black,
             pdfproducer={pdfLaTeX},
             pdfpagemode=None,
             bookmarksopen=true
             bookmarksnumbered=true]{hyperref}
\usepackage{graphicx}

\usepackage{tikz}
\usetikzlibrary{decorations.markings,arrows,decorations.pathreplacing}
\tikzset{vertex/.style={inner sep=2pt,draw}}
\tikzset{dot/.style={draw,fill,circle,inner sep=1.5pt,black}}

\allowdisplaybreaks

\newcommand\arXiv[1]{\href{http://arxiv.org/abs/#1}{\nolinkurl{arXiv:#1}}}
\newcommand\MRnumber[1]{\href{http://www.ams.org/mathscinet-getitem?mr=#1}{\nolinkurl{MR#1}}}
\newcommand\DOI[1]{\href{http://dx.doi.org/#1}{\nolinkurl{DOI:#1}}}
\newcommand\MAILTO[1]{\href{mailto:#1}{\nolinkurl{#1}}}

\usepackage{etoolbox}
\patchcmd{\section}{\scshape}{\bfseries}{}{}
\makeatletter
\renewcommand{\@secnumfont}{\bfseries}
\makeatother

\usepackage{dsfont}
\renewcommand\mathbb\mathds

\newcommand\bD{\mathbb D}

\newcommand\bK{\mathbb K}

\newcommand\bN{\mathbb N}

\newcommand\bR{\mathbb R}

\newcommand\cF{\mathcal F}

\newcommand\cO{\mathcal O}

\usepackage[mathscr]{euscript}

\DeclareMathOperator\homology{H}
\renewcommand\H{\homology}

\renewcommand\d{\mathrm d}

\newcommand\longto\longrightarrow
\newcommand\mono\hookrightarrow
\newcommand\epi\twoheadrightarrow
\newcommand\isom{\overset\sim\to}
\newcommand\<\langle
\renewcommand\>\rangle
\newcommand\sminus\smallsetminus

\newcommand\op{\mathrm{op}}
\newcommand\id{\mathrm{id}}

\DeclareMathOperator\cone{cone}
\DeclareMathOperator\End{End}
\DeclareMathOperator\Spec{Spec}
\DeclareMathOperator\maps{maps}
\DeclareMathOperator\supp{supp}
\DeclareMathOperator\Qloc{Qloc}
\DeclareMathOperator\diag{diag}

\newcommand\Frob{\mathrm{Frob}}
\newcommand\Arr{{\mathrm{Arr}}}
\newcommand\Com{\mathrm{Com}}
\newcommand\Pois{\mathrm{Pois}}
\newcommand\LB{\mathrm{LieBi}}
\newcommand\hFrob{\mathrm{hFrob}}

\newcommand\oco{{\mathrm{oco}}}
\newcommand\str{{\mathrm{str}}}

\newcommand\define[1]{\emph{#1}}
\newcommand\cat[1]{\textsc{#1}}

\title[Exact triangles, Koszul duality,  and coisotopic boundary conditions]{Exact triangles, Koszul duality,  \\ and coisotopic boundary conditions}
\author{Theo Johnson-Freyd}
\date{\today}
\email{\MAILTO{theojf@pitp.ca}}
\address{Perimeter Institute for Theoretical Physics, Waterloo, Ontario}

\begin{document}
\begin{abstract}
We develop a theory of ``arrowed'' (operads and) dioperads, which are to exact triangles as dioperads are to vector spaces.  
A central example to this paper is the arrowed operad controlling ``derived ideals'' for any operad.
The Koszul duality theory of arrowed dioperads interacts well with rotation of exact triangles, and in particular with ``exact Stars of David,'' which are pairs of exact triangles drawn on top of each other in an interesting way.  Using this framework, we give a cochain-level lift of the ``relative Poincar\'e duality'' enjoyed by oriented manifolds with boundary; moreover, our cochain-level lift satisfies a natural locality-type condition, and is uniquely determined by this property.  We discuss the meaning of the words ``relative orientation'' and ``coisotropic.''  We extend the  AKSZ construction to bulk-boundary settings with Poisson bulk fields and coisotropic boundary conditions.
\end{abstract}
\maketitle

\setcounter{subsection}{-1}

\subsection{}  The motivation for this paper came from an attempt to extend the author's ``Poisson AKSZ construction'' from \cite{PoissonAKSZ} to manifolds with boundary, where the bulk fields should be valued in a Poisson infinitesimal manifold and the boundary fields should be in a coisotropic submanifold thereof.  This requires understanding  what ``coisotropic submanifold'' should mean in infinitesimal derived geometry.  We will get to that question eventually, but most of the paper concerns abstract nonsense of operads, dioperads, and exact triangles that is likely of independent interest.  The ``Poisson'' portion of the story begins in section~\ref{derived geometry}.

\subsection{} \label{section rotation}

Let us say that an \define{arrow} is a pair of cochain complexes $A$ and $B$ together with a cochain map $f: A \to B$.  (A field $\bK$ of characteristic zero is fixed at the outset, and all cochain complexes, etc., are understood as being over $\bK$.)  It is a commonplace that the (derived or $\infty$-) category $\cat{Arrows}$ of arrows is equivalent to the  category of exact triangles (which are not truly triangles, but helices spiraling in the ``degree'' direction).  This latter incarnation provides a cube root of the \define{shift} functor $(A \overset f \to B) \longmapsto (A[1] \overset {f[1]} \to B[1])$, namely the rotation of triangles:
\begin{equation} \label{rotation}
\begin{tikzpicture}[baseline=(middle)]
  \path (150:2) node (A) {$A$} +(0,-24pt) node[black!70] (A1) {$A[1]$} (30:2) node (B) {$B$} +(0,-24pt) node[->,black!40] (B1) {$B[1]$} (-90:2) node (C1) {$C$} +(0,24pt) node (C) {$C[-1]$} +(0,-24pt) node[black!10] (C2) {$\dots$} (0,-20pt) coordinate (middle);
  \draw[->,black!10] (B1) --  (C2);
  \draw[->,black!40] (A1) -- node[auto] {$\scriptstyle f[1]$} (B1);
  \draw[->,black!70] (C1) -- node[auto] {$\scriptstyle h[1]$} (A1);
  \draw[->,thick] (B) -- node[auto] {$\scriptstyle g$} (C1);
  \draw[->,thick] (A) -- node[auto] {$\scriptstyle f$} (B);
  \draw[->,thick] (C) -- node[auto] {$\scriptstyle h$} (A);
\end{tikzpicture}
\quad\longmapsto\quad
\begin{tikzpicture}[baseline=(middle)]
  \path (150:2) node (A) {$B$} +(0,-24pt) node[black!70] (A1) {$B[1]$} (30:2) node (B) {$C$} +(0,-24pt) node[->,black!40] (B1) {$C[1]$} (-90:2) node (C1) {$A[1]$} +(0,24pt) node (C) {$A$} +(0,-24pt) node[black!10] (C2) {$\dots$} (0,-20pt) coordinate (middle);
  \draw[->,black!10] (B1) --  (C2);
  \draw[->,black!40] (A1) -- node[auto] {$\scriptstyle g[1]$} (B1);
  \draw[->,black!70] (C1) -- node[auto] {$\scriptstyle f[1]$} (A1);
  \draw[->,thick] (B) -- node[auto] {$\scriptstyle h[1]$} (C1);
  \draw[->,thick] (A) -- node[auto] {$\scriptstyle g$} (B);
  \draw[->,thick] (C) -- node[auto] {$\scriptstyle f$} (A);
\end{tikzpicture}
\end{equation}

Suppose that $A\overset f \to B$ and $X \overset p \to Y$ are arrows.  One can form a new arrow by tensoring these two together componentwise: \begin{equation} \label{tensor product of arrows} (A \overset f \to B) \quad \otimes \quad (X \overset p \to Y) \quad = \quad (A \otimes X) \overset {f \otimes p} \longto (B \otimes Y) \end{equation}
This operation of tensor product makes the category of arrows into a symmetric monoidal category, and satisfies a natural compatibility with the shift functor.  It is not, however, particularly compatible with rotation.  That said, there is a natural ``tensor multiplication'' that is rotation-invariant.  Recall that the \define{cone} of an arrow is (up to a degree convention) the value of the third vertex of the corresponding exact triangle: $$\cone\bigl(A\overset f \to B\bigr) = \Bigl(A[1] \oplus B, \, \partial = \bigl(\begin{smallmatrix} \partial_A & 0 \\ f & \partial_B \end{smallmatrix}\bigr) \Bigr).$$
Suppose that $C[-1] \overset h \to A \overset f \to B \overset g \to C$ and $Z[-1] \overset r \to X \overset p \to Y \overset q \to Z$ are two exact triangles.  Then it is not hard to see (axiom TC3 from \cite{MR1867203}) that there are canonical quasiisomorphisms
\begin{align} \label{rotation and tensoring} \cone\bigl( (A \otimes X) \overset {f \otimes p} \longto (B \otimes Y) \bigr) &\quad\simeq\quad \cone\bigl( (B \otimes Z[-1]) \overset {g \otimes r} \longto (C \otimes X) \bigr)  \\ \notag &\quad\simeq\quad \cone\bigl( (C[-1] \otimes Y) \overset {h \otimes q} \longto (A \otimes Z) \bigr).\end{align}
Indeed, up to quasiisomorphism one can substitute $C \simeq \cone(A\overset f \to B)$ and $Z \simeq \cone(X\overset p \to Y)$ and then use that there is a deformation retraction
\begin{equation*} \cone\bigl( (L \oplus M) \overset{( \begin{smallmatrix} a & 0 \\ \id & b \end{smallmatrix})}\longto (M\oplus N) \bigr) \quad \simeq \quad \cone \bigl( L \overset{b\circ a}\longto N \bigr). \end{equation*}

One can encode equation (\ref{rotation and tensoring}) in terms of an ``exact Star of David'' formed from two interlocking exact triangles:
\begin{equation} \label{star of david}
\begin{tikzpicture}[baseline=(middle),thick]
  \path 
  (150:3) node (A) {$A$}  (30:3) node (B) {$B$} (-90:3) node (C1) {$C[-1]$} ++(0,18pt) node (C) {$C$} 
  (-150:3) node (X) {$X$}  (-30:3) node (Y) {$Y$} (90:3) node (Z1) {$Z[-1]$}  ++(0,-18pt) node (Z) {$Z$}
  (0,0) coordinate (middle);
  \draw[->] (C1) -- node[auto,swap] {$\scriptstyle h$} (A);
  \draw[->] (B) -- node[auto,swap,pos=.6] {$\scriptstyle g$} (C);
  \draw[->] (A) -- node[auto,swap] {$\scriptstyle f$} (B);
  \draw[->] (Z1) -- node[auto] {$\scriptstyle q$} (X);
  \draw[->] (Y) -- node[auto,pos=.6] {$\scriptstyle r$} (Z);
  \draw[->] (X) -- node[auto] {$\scriptstyle p$} (Y);
\end{tikzpicture}
\end{equation}
The rotation-invariant cochain complex from equation (\ref{rotation and tensoring}) is the cone of the tensor product of any pair of parallel edges.  Note that the two triangles in (\ref{star of david}) are oriented in different directions: thus, ``forward'' rotation of the clockwise $ABC$ triangle, as in~(\ref{rotation}), must be complemented by ``backward'' rotation of the counterclockwise $XYZ$ triangle.  For example, rotating each triangle three times in complementary directions gives the obvious isomorphism
$$ \cone\bigl( (A \otimes X) \overset {f \otimes p} \longto (B \otimes Y)\bigr) \quad \cong \quad \cone\bigl( (A[1] \otimes X[-1]) \overset {f[1] \otimes p[-1]} \longto (B[1] \otimes Y[-1])\bigr). $$

\subsection{} \label{section arrowed operad}

We will be interested in situations where an arrow $A \overset f \to B$ is equipped with some extra algebraic structure.  The most popular language in which to describe types of algebraic structures is the language of (colored) operads.  We will need a more flexible language that can accommodate many-to-many operations.  There are essentially two such languages: \define{dioperads} describe situations in which all compositions are ``tree-shape'' whereas \define{props} describe situations in which compositions can be parameterized by more complicated graphs.  
The free/forget adjunction between dioperads and props does not provide an adjunction of $\infty$-categories \cite{MR2320654}, and so the choice of language really does matter when trying to understand homotopy algebra.  

The examples that conclude this paper are dioperadic, and so that is the context we will focus on.  The precise definition of (uncolored) {dioperad} is in \cite{MR1960128}; roughly speaking, a (dg) dioperad $P$ consists of cochain complexes $P(k;k')$ of ``operations with $k$ inputs and $k'$ outputs'' for each pair $(k,k') \in \bN^2$, which should carry commuting actions of the symmetric groups $S_k$ and $S_{k'}$, and ``compositions'' parameterized by directed trees. 

Any associative algebra $A$ gives an example of a dioperad by setting $P(1;1) = A$ and $P(k;k')= 0 $ if $(k,k') \neq (1,1)$.  Let $\Arr$ denote the algebra whose modules are arrows $A \overset f \to B$: the standard basis for $\Arr$ is $\{\pi_A,\pi_B,f\}$ and the only non-zero multiplications in this basis are $\pi_A^2 = \pi_A$, $\pi_B^2 = \pi_B$, and $f\pi_A = \pi_B f = f$.  The main algebraic objects of interest in this paper are \define{arrowed dioperads}, which are dioperads $P$ equipped with a homomorphism $\Arr \to P$.  Dioperads parameterize algebraic structures that can be carried by cochain complexes: if $V$ is a cochain complex, there is a dioperad $\End(V)$ with $\End(V)(k;k') = \hom(V^{\otimes k},V^{\otimes k'})$, and a \define{$P$-algebra structure on $V$} is a homomorphism $P \to \End(V)$.  Arrowed dioperads, in turn, parameterize algebraic structures that can be carried by arrows: given an arrow $A \overset f\to B$, let $\End\bigl(A \overset f \to B\bigr)(k;k') = \hom\bigl((A\oplus B)^{\otimes k},(A\oplus B)^{\otimes k'}\bigr)$ with the obvious structure as an arrowed dioperad; then an \define{arrowed algebra} is a homorphism $P \to \End(A\overset f \to B)$ of arrowed dioperads.

We will generally work with arrowed dioperads not in terms of dioperad homomorphisms $\Arr \to P$ but in an equivalent somewhat-unpacked form.  Any arrowed dioperad is equivalent to the following ``bicolored'' data:
\begin{itemize}
  \item two colors --- generally ``A'' and ``B'';
  \item for each four-tuple $(m,n,m',n')$, a cochain complex $P(m,n;m',n')$ of operations with $m$ A-colored inputs, $n$ B-colored inputs, $m'$ A-colored outputs, and $n'$ B-colored outputs;
  \item a distinguished \define{arrowing} $f\in P(1,0;0,1)$;
  \item an action of $S_m \times S_{m'} \times S_n \times S_{n'}$ on $P(m,n;m',n')$;
  \item compositions parameterized by directed trees whose edges are each colored~``A'' or~``B.''
\end{itemize}
For example, $\End\bigl(A\overset f \to B\bigr)(m,n;m',n') = \hom(A^{\otimes m} \otimes B^{\otimes n},A^{\otimes m'} \otimes B^{\otimes n'})$.
Non-arrowed dioperads will variously be called ``uncolored,'' ``singly-colored,'' or ``plain.''

It will sometimes be convenient to assume our arrowed operad satisfies various technical conditions. 
A dioperad (arrowed or plain) is \define{locally finite} if it is finite-dimensional in each arity.  
A plain dioperad $O$ is \define{open and coopen}, abbreviated \define{oco}, if $O(k;0) = O(0;k) = 0$ and $O(1;1)$ is one-dimensional spanned by the identity element.  (The name comes from the case of ``open Frobenius algebras''; see section~\ref{bv tensor product}.)  
Similarly, we will say that an arrowed dioperad is \define{oco} if it has no operations of arity $(0,0;m',n')$ or $(m,n;0,0)$ and the only operations of total arity $\leq 2$ are the linear multiples of the two identity elements and the arrowing (which should be non-zero).  
If $P$ is oco, we let $\underline P$ denote the collection of all operations in $P$ with total arity at least $3$. 
The data of the arrowing $f$ is retained by the collection of maps $f\circ : \underline P(m,n;m',n') \to \underline P(m,n;m'-1,n'+1)$ and $\circ f:\underline P(m,n;m',n') \to \underline P(m+1,n-1;m',n')$.

Suppose that $G = \{G(m,n;m',n')\}_{m,m',n,n'\in \bN}$ is a collection of $S_m \times S_{m'} \times S_n \times S_{n'}$-modules; elements of $G$ can be visualized as bicolored corollas.  Then $G$ generates a free arrowed dioperad $\cF(G)$.  This is also the free bicolored dioperad generated by $G \sqcup \{f\}$.  A generating set $G$ is \define{locally finite} and \define{oco} if all $G(m,n;m',n')$ are finite-dimensional and $G$ contains neither generators of total arity $\leq 2$ nor generators with $(m+n)(m'+n') = 0$.  If $G$ is locally finite and oco, then so is $\cF(G)$.

Finally, an arrowed dioperad $P$ is \define{quasifree} if it is free upon forgetting the differential, and moreover the generating set $G$ can be well-ordered so that for each generator $x\in G$, $\partial x \in P$ is in the sub arrowed dioperad generated by generators strictly before $x$ for the well-ordering.  
The existence of a well-ordering is automatic when $G$ is oco, but in general it is very strong. It implies, by the usual arguments, that quasifree arrowed dioperads satisfy the left-lifting property against surjective quasiisomorphisms. Every oco arrowed dioperad has a quasifree resolution, for example its double bar dual as described in section~\ref{bar dual}.
More general model-categorical control can be provided by following the construction from the appendix of \cite{MR2572248}.

\subsection{Example}   \label{Oto}

A dioperad $O$ is an \define{operad} if $O(k;k') = 0$ unless $k'=1$. An \define{arrowed operad} is an arrowed dioperad such that $P(m,n;m',n') = 0$ unless $m'+n' = 1$.   Suppose $O$ is an (uncolored) operad.   Then there is an arrowed operad $O^{\to,\str}$ whose algebras are pairs of $O$-algebras connected by a (strict) homomorphism.  It satisfies 
$$ O^{\to,\str}(m,n;m',n') \cong \begin{cases} 
  O(m;1), & (n,m',n') = (0,1,0) \\
  O(m+n;1), & (m',n') = (0,1) \\
  0, & \text{else}
 \end{cases} $$
 with the obvious compositions.  The arrowing is the element of $O^{\to,\str}(1,0;0,1)$ corresponding to the identity element in $O(1;1)$.  If $O$ is oco, so is $O^{\to,\str}$.  Any quasiisomorphism $O \isom O'$ induces a quasiisomorphism $O^{\to,\str} \isom (O')^{\to,\str}$.  Define an \define{$\infty$-morphism} of $O$-algebras to be an arrowed algebra for any quasifree arrowed operad $O^\to$ resolving $O^{\to,\str}$.

In order to make contact with the standard meaning of ``$\infty$-morphism'' described in \cite[Section~10.2]{MR2954392}, suppose that $O$ is \define{quadratically quasifree} in the sense that the derivative of each generator $x \in G$ is quadratic in generators, and suppose furthermore that the complex $G$ of generators is non-zero only in arities $(k,1)$ for $k\geq 2$. (We will discuss bar duality in section~\ref{bar dual}; up to some set-theoretic issues that can be corrected by working also with co-structures, quadratically quasifree operads are precisely the bar duals of operads.) 
An `$\infty$-morphism of $O$-algebras'' $A$ and $B$ in the sense of \cite[Section~10.2]{MR2954392} then unpacks to a linear map $f : A \to B$ and, for each generator $x\in G(k)$, a map $x^\to : A^{\otimes k} \to B$ of cohomological degree $\deg(x^\to) = \deg(x) - 1$, subject to a cohomological constraint. Specifically, letting $x^A: A^{\otimes k} \to A$ and $x^B: B^{\otimes k} \to B$ denote the values of the generators $x \in G$ in the algebras $A$ and $B$, the constraint has the form
$$ \partial_B \circ x^\to \pm x^\to \circ \partial_{A^{\otimes k}} = f \circ x_A - x_B \circ f^{\otimes k} + \text{composition of generators of arity $<k$,}$$
where the sign is determined by the Koszul sign rule.

This ``Loday--Vallette'' description of an $\infty$-morphism can be packaged into a quasifree arrowed operad $O^{\to,LV}$ with generators
$$ G^{\to, LV}(m,n;m',n') = \begin{cases}
  G(k;1), & (m,n,m',n') = (k,0,1,0) \\
  G(k;1), & (m,n,m',n') = (0,k,0,1) \\
  G(k;1)[1], & (m,n,m',n') = (k,0,0,1) \\
  0, & \text{else}
\end{cases}$$
By construction, there is a surjection $O^{\to,LV} \to O^{\to,\str}$ acting as the identity on the generators of arities $(k,0,1,0)$ and $(0,k,0,1)$ and annihilating the generators in arity  $(k,0,0,1)$. (One could say that it sends those generators to the relation $f \circ x_A = x_B \circ f^{\otimes x}$.)
This surjection is manifestly an isomorphism in all arities except $(m,n;0,1)$.
To see that it is a quasiisomorphism in arity $(m,n;0,1)$, one can choose a filtration so that the derivatives read
$$ [\partial,x^\to] = f \circ x^A + \text{lower order}.$$
After ``canceling'' generators $x^\to$ and $x^A$, one is left just with diagrams of the form $($composition of $x^B$s$) \circ f^{\otimes \dots}$, which are in bijection with elements of $O^{\to,\str}(m,n;0,1)$ as desired.

Since $O^{\to,LV} \to O^{\to,\str}$ is a quasifree resolution, $O^{\to,LV}$ provides provides a model of $\infty$-morphisms in our sense. Furthermore, standard model-categorical nonsense guarantees that any other model is equivalent: for any choice of $O^\to$, each $O^\to$-algebra corresponds to a homotopically-unique $O^{\to,LV}$-algebra.

The construction of $O^{\to,LV}$ and of the quasifree resolution $O^{\to,LV} \to O^{\to,\str}$ is outlined in~\cite{MarklTalk}, where the example $O = A_\infty$ is worked out in detail.

\subsection{}\label{O ideal}

We will use heavily the following construction.  Let $O$ be a (singly-colored) dioperad.  Define the arrowed dioperad $O^{\triangleleft,\str}$ by
\begin{equation}\label{strict ideal} O^{\triangleleft,\str}(m,n;m',n') = \begin{cases} 0, & m=n'=0 \\ O(m+n;m'+n'), & \text{else} \end{cases}\end{equation}
Composition uses the composition in $O$.  To see that composition never results in the ``$m=n'=0$'' case,  replace every A-colored edge by a ``forward'' arrow and  every B-colored edge by a ``backward'' arrow; then the rule for $O^{\triangleleft,\str}$ is that no vertex should be a source.
  A directed tree with no sources necessarily has an incoming leaf, proving the claim.  
  (If we were using props, which allow graph-based compositions, then $O^{\triangleleft,\str}(m,n;m',n')$ as defined would not be closed under composition.)
The arrowing $f\in O^{\triangleleft,\str}(1,0;0,1)$ is the element corresponding to $\id \in O(1;1)$.
Clearly $(-)^{\triangleleft,\str}$ preserves quasiisomorphisms.  We will refer to $O^{\triangleleft,\str}$ as \define{strict $O$-(co)ideals}.

Suppose now that $O$ is quasifree with generators $G$.  Let $O^{\triangleleft}$ denote the quasifree arrowed dioperad with generators
$$ G^{\triangleleft}(m,n;m',n') = \begin{cases} 0, & m=n'=0 \\ G(m+n;m'+n')[m+n'-1], & \text{else} \end{cases}$$
and differential defined as follows.  Given a directed tree $\Gamma \in O(m+n;m'+n')$ with vertices labeled by generators $x \in G$, let $\Gamma_{m,n;m',n'}$ denote the sum of bicolored trees formed from $\Gamma$ by coloring the first $m$ incoming and first $m'$ outgoing edges A and the last $n$ incoming and last $n'$ outgoing edges B, and summing over all ways to color interior edges A or B such that no vertex has all incoming edges colored B and all outgoing edges colored A.  In particular, if $m=n'=0$, the sum is empty, and so $\Gamma_{0,m';n,0} = 0$.  When $\Gamma = x\in G(m+n;m'+n')$ is a generator of $O$,  $x_{m,n;m',n'}$ is the corresponding generator of $O^{\triangleleft}$ with the given arity.  The differential on $O^{\triangleleft}$ is defined by
\begin{equation}\label{O ideal eqn} \partial(x_{m,n;m',n'}) = (\partial x)_{m,n;m',n'} + \sum_{j=1}^m (-1)^j x_{m-1,n+1;m',n'} \circ_{j} f - \sum_{j'=1}^{n'} (-1)^{j'} f \circ_{j'} x_{m,n;m'+1,n'-1}. \end{equation}
Algebras for $O^\triangleleft$ are \define{$O$-$\infty$-(co)ideals}; Theorem~\ref{O ideal presents O strict ideal} shows that $O^\triangleleft$ presents the ``$\infty$-algebras'' for $O^{\triangleleft,\str}$.

We now begin to justify the name ``(co)ideal''; complete justification is given by Theorem~\ref{O ideal justification}.  Suppose that $O$ is a quasifree operad.  Then $O^\triangleleft$ has the following interpretation.  First, $O^\triangleleft(0,k;0,1) = O(k;1)$, and so the generators $x_{0,k;0,1}$ make the ``B'' component of any $O^\triangleleft$-algebra into an $O$-algebra.  Next, the generators $x_{1,k-1;1,0}$ make $A$ into an $\infty$-$B$-module.  
The generators $x_{1,k-1;0,1}$ make $f: A \to B$ into an $\infty$-morphism of $A$-modules.  
 Given two elements $a_1,a_2 \in A$ and some elements $b_3,\dots,b_k \in B$, one can now consider two different ``multiplications'' valued in $A$: act on $a_1$ by $f(a_2), b_3,\dots,b_k$, or act on $a_2$ by $f(a_1), b_3,\dots,b_k$. The generators $x_{2,k-2;1,0}$ identify these two $A$-valued ``multiplications,'' and the generators $x_{2,k-2;0,1}$ then say that $f$ is compatible with this identification.
 And so on.  
    The end result is that $O^\triangleleft$ presents pairs $(A,B)$ where $B$ is an $O$-algebra and $A \triangleleft B$ is an ideal in some homotopical sense.
As for the ``co,'' note that the construction of $O^\triangleleft$ is manifestly symmetric under taking a dioperad $O$ to its opposite $O^{\mathrm{op}}(k;k') = O(k';k)$.  The word ``co-operad'' being already taken, let us say that an \define{opposite operad} (an ``erad''?)\ is a dioperad which vanishes except in arity $(1;k')$.  If $O$ is an opposite operad, then the same logic implies that $O^\triangleleft$ parameterizes some homotopical version of $O$-coideals.

\subsection{Theorem} \label{O ideal presents O strict ideal} \emph{There is a surjective quasiisomorphism $O^\triangleleft \overset\sim\epi O^{\triangleleft,\str}$ sending $x_{0,n;m',1}$ and $x_{1,n;m',0}$ to the corresponding generators $x_{n;m'+1}$ and $x_{n+1;m'}$ of $O^{\triangleleft,\str}$ and annihilating all other generators.}

\subsection*{Proof}
We first must check that the claimed map $O^\triangleleft \to O^{\triangleleft,\str}$ exists; it is then manifestly surjective.  For existence, it suffices to check that for any generator $x_{m,n;m',n'}$ annihilated by the map, $\partial (x_{m,n;m',n'})$ is also annihilated.  The generators annihilated by the map are those with $m+n'\geq 2$.  If $m = n' = 1$, the term  in equation (\ref{O ideal eqn}) including $f$s is
$$ -x_{0,n+1;m',1} \circ f + f\circ x_{1,n;m'+1,0} \quad \mapsto \quad -x_{n+1;m'+1} + x_{n+1;m'+1}=0.$$
If $(m,n') = (2,0)$, then we have
$$ -x_{1,n+1;m',0} \circ_1 f + x_{1,n+1;m',0} \circ_2 f \mapsto -x_{n+2;m'} + x_{n+2;m'}=0,$$
and similarly for $(m,n') = (0,2)$.  
If $m+n' \geq 3$ then the term including $f$s is manifestly annihilated.  Finally, for any $x$ every summand in the term $(\partial x)_{m,n;m',n'}$ is annihilated.  Indeed, as in section~\ref{O ideal}, replace every A-colored edge by a ``forward'' arrow and  every B-colored edge by a ``backward'' arrow.  Then each summand in $(\partial x)_{m,n;m',n'}$ is a directed tree with $m+n'$ outgoing leaves and the rest incoming and no sources.  In a directed tree with no sources, if there is exactly one outgoing leaf, then every vertex has exactly one outgoing edge; if there are at least two outgoing leaves, then there is at least one vertex with at least two outgoing edges.  This verifies that the  map $O^\triangleleft \to O^{\triangleleft,\str}$ is well-defined.

To check that it is a quasiisomorphism, it suffices to consider the case when all $\partial x$s vanish: they are of subleading order in equation~(\ref{O ideal eqn}), and turning on subleading  terms can never spoil a quasiisomorphism.  
We therefore assume for the remainder of the proof that $O$ is free and not just quasifree.

A general element of $O^\triangleleft(m,n;m',n')$ is a linear combination of bicolored directed trees with vertices labeled by generators of $O$ and the occasional bivalent vertex labeled ``$f$.''  Any such tree determines an element of $O$ by ignoring the ``$f$'' vertices. The element of $O$ and the corresponding element of $O^\triangleleft$ differ in degree by $1-(m+n')+\text{the number of ``$f$'' vertices}$.  For comparison, a typical element of $O^{\triangleleft,\str}$ looks like a tree in the generators of $O$ with only its leaves colored.  Thus we can think of a tree $\Gamma \in O^\triangleleft$ as a tree $[\Gamma] \in O^{\triangleleft,\str}$ together with a coloring of its edges by the colors ``A,'' ``B,'' and ``$f$,'' subject to two rules: (1) B-colored incoming leaves on $[\Gamma]$ are B-colored on $\Gamma$ and A-colored incoming leaves on $[\Gamma]$ are either A- or $f$-colored on $\Gamma$, and A-colored outgoing leaves on $[\Gamma]$ are A-colored on $\Gamma$ and B-colored outcoming leaves on $[\Gamma]$ are either B- or $f$-colored on $\Gamma$; (2) no vertex has all B- and $f$-colored incoming edges and all A- and $f$-colored outgoing edges.  

The differential does not change $[\Gamma]$.  Indeed, equation~(\ref{O ideal eqn}) can be re-read as saying that $\partial$ acts on colored edges, not on vertices, where $\partial(\text{A-colored edge}) = \text{$f$-colored edge} = -\partial(\text{B-colored edge})$ (and $\partial(\text{$f$-colored edge}) = 0$), with the caveat that the derivative vanishes if the coloring violates the rules.  Thus it suffices to prove that the complex of all colorings of a fixed $[\Gamma]$ has just one-dimensional cohomology.

Fix $[\Gamma]$ and suppose that its edges have been partially colored compatibly with the rules in the sense that there exists a compatible full coloring.  Let $e$ be a not-yet-colored interior edge.  If ``$f$'' is an allowed coloring for $e$, then both ``A'' and ``B'' are allowed colorings.  Otherwise, exactly one of ``A'' or ``B'' is allowed.  In either case, the complex of allowed colorings of $e$ has one-dimensional cohomology (in the degree without an ``$f$'' edge).
To compute the cohomology, then, we can ignore all interior colorings.  Rule (1)  dictates the colorings of some of the exterior edges, but allows others two choices, provided they can be made compatibly with rule (2).  Let $L$ denote the set of A-colored incoming exterior edges and B-colored outgoing exterior edges in $[\Gamma]$, so that $|L| = m+n'$.  Rule (2) allows for all but one of these to be colored ``$f$'' but not all of them; otherwise it provides no rules, after the complex of colorings of interior edges has been contracted to its homology.  Thus the complex of allowed colorings of $L$ looks like the totalization of a cube $(\bK \overset \partial  \to \bK)^{\otimes L}$ without its ``all $f$s'' vertex.  Its homology is therefore one-dimensional supported in the degree for which there are $|L|-1$ many $f$s. \qed

\subsection{}\label{rotation of arrowed operads} 

Let $P$ be an arrowed dioperad and $\cat{Alg}_P$ the category of arrowed $P$-algebras.  There is a forgetful functor $\operatorname{Forget}: \cat{Alg}_P \to \cat{Arrows}$, and the philosophy of the Barr--Beck theorem says that $P$ can be reconstructed from $\operatorname{Forget}$.  Consider composing $\operatorname{Forget}$ with an autoequivalence of $\cat{Arrows}$.  One produces a new functor $\operatorname{Forget}': \cat{Alg}_P \to \cat{Arrows}$ with the same formal properties as $\operatorname{Forget}$, and so the Barr--Beck philosophy suggests that one can reconstructs some new arrowed dioperad $P'$.  For example, composing $\operatorname{Forget}$ with the shift functor $(A\overset f \to B) \mapsto (A[1] \overset {f[1]}\to B[1])$ corresponds to \define{shearing} $P$ to the arrowed dioperad $P\langle 1\rangle$ defined by $P\langle 1\rangle (m,n;m',n') = P(m,n;m',n')[m'+n'-m-n]$.  (Here and throughout we adopt the convention that if $V$ is an $S_k$-module, then $V[k]$ is given the $S_k$-module structure coming from the diagonal action on $V \otimes (\bK[1])^{\otimes k}$, where $S_k$ acts on $(\bK[1])^{\otimes k}$ by permuting \emph{with the Koszul sign rule} the $k$ copies of $\bK[1]$.)

As mentioned already in section \ref{section rotation}, $\cat{Arrows}$ admits much more interesting autoequivalences: $(A \overset f \to B) \mapsto (B \to \cone(f))$ and its inverse $(A \overset f \to B) \mapsto (\cone(f)[-1] \to A)$.
We now work out how these operations act on quasifree arrowed dioperads. 
 Specifically, we look for quasifree arrowed dioperads $\Theta P$ and $\Theta^{-1}P$ such that if $C[-1] \overset h \to A \overset f \to B \overset g \to C$ is an exact triangle, then $P$-algebra structures on $A \overset f \to B$ are the same as $\Theta P$-algebra structures on $B \overset g\to C$ and as $\Theta^{-1}P$-algebra structures on $C[-1] \overset h \to A$.

To preserve the notation, let us call $\Theta P$'s colors ``B'' and ``C'' in that order, and the arrowing~``$g$.''
The strategy is to take the generators of $P$ and substitute $A = \cone(B \overset g \to C)[-1] = (B \oplus C[-1], \partial = \partial_B + \partial_C + g)$.
 Let $x$ be a generator of $P$ with arity $(m,n;m',n')$.  
 Choose  subsets $L \subseteq \{1,\dots,m\}$ and $L' \subseteq \{1,\dots,m'\}$ and set $\ell = |L|$ and $\ell' = |L'|$.  
 Then $\Theta P$ will have a generator $x_{L;L'}$ corresponding to $x$ with the A-colored inputs in $L$ now colored $C$,
  the remaining $m-\ell$ A-colored inputs now colored B, 
  the A-colored outputs in $L'$ now colored C,
  and the remaining $m'-\ell'$ A-colored outputs now colored B.
 The degree of $x_{L;L'}$ shifts from that of $x$ by $\ell' - \ell$. 
 By pre- and post-composing $x$ with permutations of the inputs and outputs, we can assume that $L = \{1,\dots,\ell\}$ and $L' = \{1,\dots,\ell'\}$;
 thus if we adopt the convention that the generators are closed under the symmetric group actions, then we can simply write $x_{\ell;\ell'}$, $\ell = 0,\dots,m$, $\ell' = 0,\dots,m'$, for the new generators, rather than $x_{L;L'}$. 
  More precisely, consider the complex $G(m,n;m',n')$ of generators of $P$ with arity $(m,n;m',n')$.  This complex carries an action of the product of symmetric groups $S_m \times S_n \times S_{m'} \times S_{n'}$. 
   For each $\ell \leq m$ and $\ell'\leq m'$, we restrict to $S_\ell \times S_{m-\ell}  \times S_n \times S_{\ell'} \times S_{m'-\ell'} \times S_{n'} $
  and then induce to $S_{m + n - \ell} \times S_\ell \times S_{m'+n'-\ell'}\times S_{\ell'} $.  This gives the new complex of $x_{L;L'}$s for $(|L|,|L'|) = (\ell,\ell')$.   

The data of $P$ consists of the generators $x$ together with their derivatives $\partial x$, each of which is a sum of directed trees in $f$ and earlier generators, with edges colored A or B. 
To describe $\Theta P$, we need simply to present a formula for $\partial(x_{\ell;\ell'})$.  Suppose $\Gamma \in P(m,n;m',n')$ is a bicolored directed tree.  Let $\Gamma_{\ell;\ell'}$ denote the sum of trees formed as follows.  
The first $\ell$ A-colored incoming and first $\ell'$ A-colored outgoing leaves switch color to C, and the remainder switch to B.
The sum ranges over all ways to convert interior A-colored edges to either B or C.  
In each summand, after recoloring one gets a tree with edges colored either B or C and with vertices coming from generators $x$ of $P$, and with the occasional ``$f$'' vertex.  
The ``$x$'' vertices are switched to the corresponding generator of $\Theta P$ according to the local coloring.  Before recoloring, the ``$f$'' vertices had A-colored input and B-colored output; after recoloring, they have 
either B- or C-colored input and B-colored output.
If the recolored ``$f$'' vertex has B-colored input,  simply smooth out the vertex, replacing it by the identity B-colored edge.
If the recolored ``$f$'' vertex has C-colored input, then replace it by $0$, thereby annihilating the whole summand.  
Then
$$ \partial(x_{\ell;\ell'}) = (\partial x)_{\ell;\ell'} + \sum_{j=1}^{m-\ell} x_{\ell+1;\ell'} \circ_j g = \sum_{j'=1}^{\ell'} g \circ_{j'} x_{\ell;\ell'-1}$$
where $\circ_j$ indicates that the composition occurs at the $j$th leaf.

A similar description presents $\Theta^{-1}P$.  Indeed, arrowed dioperads admit a manifest symmetry reversing incoming and outgoing directions (and reversing colors): the \define{opposite} of an arrowed dioperad $P$ is $P^{\mathrm{op}}(m,n;m',n') = P(n',m';n,m)$.  Then $\Theta^{-1}P = (\Theta (P^{\mathrm{op}}))^{\mathrm{op}}$.

It is clear from the definition that if $P$ is oco, so are $\Theta P$ and $\Theta^{-1} P$.

\subsection{Theorem} \label{O ideal justification}

\emph{If $O$ is a quasifree operad, there is a surjective quasiisomorphism $\Theta O^\triangleleft \overset\sim\epi O^{\to,\str}$.}
\\ \mbox{}

Equivalently, if $O$ is an opposite operad, then strict morphisms of $O$-algebras are parameterized by an arrowed opposite operad $O^{\to,\str}$ defined analogously to Example~\ref{Oto}, and the Theorem then asserts that $\Theta^{-1}O^\triangleleft \simeq O^{\to,\str}$.  Together these justify the name ``$O$-(co)ideal'' in section~\ref{O ideal}: if $O$ is an (opposite) operad, $O^\triangleleft$-algebras are then the homotopy-(co)kernels of $O^\to$-algebras.
A related result in the fully dioperadic case is in Example~\ref{cone of an ideal}.

\subsection*{Proof}

Since $O$ is an operad, its generators are all of arity $(k;1)$ for various $k$.  A generator $x \in O(k;1)$ effects the following generators of $\Theta O^\triangleleft$:
\begin{align*}
  x_{B,m,\ell} & = (x_{m,k-m;1,0})_{\ell;0} \in \Theta O^\triangleleft(k-\ell,\ell;1,0) [1-m+\ell], \quad m=1,\dots,k,\quad \ell = 0,\dots,m, \\
  x_{C,m,\ell} & = (x_{m,k-m;1,0})_{\ell;1} \in \Theta O^\triangleleft(k-\ell,\ell;0,1) [-m+\ell], \quad m=1,\dots,k,\quad \ell = 0,\dots,m, \\
  x'_{B,m,\ell} & = (x_{m,k-m;0,1})_{\ell;0} \in \Theta O^\triangleleft(k-\ell,\ell;1,0) [-m+\ell], \quad m=0,\dots,k,\quad \ell = 0,\dots,m. 
\end{align*}
The differentials are
\begin{gather*}
  \partial x_{B,m,\ell} = 
  \underbrace{(\#)  x_{B,m-1,\ell}}_{\text{if }m\geq 2} +  \sum_i  x_{B,m,\ell+1} \circ_i g + (\dots) ,\\
  \partial x_{C,m,\ell} = 
  \underbrace{(\#) x_{C,m-1,\ell}}_{\text{if }m\geq 2} +   \sum_i x_{C,m,\ell+1} \circ_i g - g \circ x_{B,m,\ell} + (\dots) ,\\
  \partial x'_{B,m,\ell} = \underbrace{x_{B,m,\ell}}_{\text{if }m\geq 1} + 
  (\#) x'_{B,m-1,\ell} +  \sum_i x'_{B,m,\ell+1} \circ_i g + (\dots) ,
\end{gather*}
where following section~\ref{rotation of arrowed operads}  we call the arrowing ``$g$.''
In each line the $(\dots)$ terms are subleading corrections coming from the differential on $O$.  The  $(\#)$ term in each line  counts (with signs) the number of ``$\circ f$s'' from the $(-)^\triangleleft$ construction that survive the $\Theta$ construction.

Call the copy of $x$ in $O^{\to,\str}(k-\ell,\ell;0,1)$ by the name $x_{C,\ell}$, and the copy of $x$ in $O^{\to,\str}(k,0;1,0)$ by the name $x_B$.  The map $\Theta O^\triangleleft \to O^{\to,\str}$ annihilates all generators except
$$ x_{B,1,0} \mapsto x_B, \quad x'_{B,0,0} \mapsto -x_B, \quad x_{C,\ell,\ell} \mapsto x_{C,\ell}, \ell=1,\dots,k. $$
Modulo subleading ``$(\dots)$'' terms, this map clearly sends $\partial x'_{B,1,0} \mapsto x_B - x_B$, $\partial x_{C,1,0} \mapsto g \circ x_B - x_{C,1} \circ g$, 
 and, for $\ell \geq 2$, $\partial x_{C,\ell,\ell-1} \mapsto x_{C,\ell} \circ g - x_{C,\ell-1}$, all of which vanish in $O^{\to,\str}$. All other generators' derivatives vanish under the map just for reasons of indexes.  Restoring the ``$(\dots)$'' terms has the effect only of allowing $\partial x_{B,1,0}$, $\partial x'_{B,0,0}$, and $\partial x_{C,\ell,\ell}$ to be the corresponding copies of $\partial x$ in the appropriate arities; the ``$(\dots)$'' terms in the derivatives of any other generator necessarily contains generators annihilated by the map.  This verifies that there is a map $\Theta O^\triangleleft \to O^{\to,\str}$, and it is obviously a surjection.

As in the proof of Theorem~\ref{O ideal presents O strict ideal}, to check that the map is a quasiisomoprhism it suffices to consider the case when $O$ is free and not just quasifree.  Then
 the differential on $\Theta O^\triangleleft$ still has a leading order term which is linear in generators and a subleading term containing ``$\circ g$'' and ``$g\circ$'' corrections.  We study the linear-in-differentials part first; with it alone, the generators by themselves form a complex.  This linear-in-differentials complex is a direct sum of two pieces, depending on the color of the output.  Since $\partial x'_{B,m,\ell} = x_{B,m,\ell} + \dots$ for $m\geq 1$, the complex of B-output generators has only one cohomology class, represented by $x'_{B,0,0}$.  

The complex of generators $x_{C,m,\ell}$ splits as a direct sum indexed by $\ell=0,\dots,k$.  The generator $x_{C,k,k}$ is alone and represents a cohomology class.  For $\ell \in \{1,\dots,k\}$, the complex is exact: each basis vector for that summand is the generator $x$ with $\ell$ of its inputs C-colored and of the remaining $k-\ell$ B-colored inputs, $m-\ell$ of them chosen as ``originally A'' and $k-m$ as ``originally B''; the differentials are just the maps changing an ``originally A'' input into an ``originally B'' input.  For $\ell=0$, the complex of $x_{C,m,0}$s would be exact if there were a copy of $x_{C,0,0}$, but that element does not appear.  It follows that the $\ell=0$ complex has one-dimensional cohomology represented by (some combination of permutations of) $x_{C,1,0}$.

Finally, turning on the ``$\circ g$'' terms, the cohomology class represented by $x_{C,1,0}$ becomes a homotopy between $x_{C,k,k} \circ g^{\otimes k}$ and $-g\circ x'_{B,0,0}$.  But since $O$ is free, the arrowed operad $O^{\to,\str}$ can be generated by the elements $x_{C,k}$ and $x_{B}$ with the only relation $x_{C,k}\circ g^{\otimes k} = g\circ x_B$, and so the map $\Theta O^\triangleleft \to O^{\to,\str}$ is a quasiisomorphism. \qed

\subsection{}  %\label{bv tensor product} 

Any arrowed dioperad $P$ immediately provides two singly-colored dioperads $\alpha P = P(-,0;-,0)$ and $\beta P = P(0,-;0,-)$ by restricting to either color.  A third dioperad is 
\begin{equation} \label{cone of an arrowed dioperad} \gamma P(k;k') = \bigoplus_{\substack{k=m+n,\\ k'=m'+n'}} \binom{k}{m} \binom{k'}{m'} P(m,n;m',n')[m'-m], \quad  \partial_{\gamma P} = \partial_P + [f,-]  ,\end{equation}
where the prefactor $\binom km \binom{k'}{m'}$ records that the $S_m \times S_{k-m} \times S_{m'} \times S_{k'-m'}$-module should be induced up to a $S_k \times S_{k'}$-module.
The letter ``$\gamma$'' stands for ``cone,'' because $P$-algebra structures on an arrow $A \overset f \to B$ induce $\gamma P$-algebra structures on $\cone(A \overset f \to B)$.  
By construction, there are equivalences $\gamma P \simeq \beta\Theta P \simeq \alpha \Theta^{-1} P\langle 1 \rangle$.

For example, if $O$ is an oco dioperad, $\alpha O^{\triangleleft,\str} \cong \beta O^{\triangleleft,\str} \cong O$.  In general, $\alpha( O^{\triangleleft,\str})$ and $\beta (O^{\triangleleft,\str})$ look like $O$ without its operations with either zero inputs or zero outputs, respectively.  If $O$ is an operad, $\alpha O^{\to,\str} \cong \beta O^{\to,\str} \cong O$.  In Example~\ref{cone of an ideal} we will compute $\gamma O^{\triangleleft,\str}$.

\subsection{}  \label{bv tensor product} 

There are various ways to combine (singly colored or arrowed) dioperads.  The most basic is the \define{Boardman--Vogt tensor product} $P \boxtimes Q$ of (singly colored) dioperads $P$ and $Q$:
$ (P \boxtimes Q)(k;k') = P(k;k') \otimes Q(k;k') $.
The defining property of $P\boxtimes Q$ is that if $V$ is a $P$-algebra and $W$ is a $Q$-algebra, then $V \otimes W$ is a $P \boxtimes Q$-algebra.  

The unit for $\boxtimes$ is the dioperad $\Frob$ of (commutative and cocommutative) \define{Frobenius algebras}, satisfying $\Frob(k;k') = \bK$ for all $k,k'$ (with the trivial symmetric group actions and the obvious composition).  The $\boxtimes$-invertible dioperads parameterize \define{shifted} Frobenius algebras.  Given integers $d,d'$, a cochain complex $V$ is a \define{$(d,d')$-shifted Frobenius algebra} if $V[d]$ is a commutative algebra and $V[d']$ is a cocommutative coalgebra, and the multiplication and comultiplication satisfy the Frobenius relation.  For example, the cohomology of a compact oriented $d'$-dimensional manifold is a $(0,d')$-shifted Frobenius algebra.
The corresponding dioperad $\Frob_{d,d'}$ satisfies
$ \Frob_{d,d'}(k;k') = \bK[dk-d'k'-d+d'] $.
Tensoring with the dioperads of shifted Frobenius algebras extends the ``shearing'' operation from section~\ref{rotation of arrowed operads}: when $d=d'$, $P \boxtimes \Frob_{d,d} = P\langle -d \rangle$.

The Boardman--Vogt tensor product of locally finite oco dioperads is again locally finite and oco.    The unit oco dioperad $\Frob^{\oco}$ parameterizes \define{nonunital and non-counital}, aka \define{open and coopen}, Frobenius algebras, and satisfies $\Frob^\oco(k;k') = \bK$ if $kk' \neq 0$ and $0$ if $kk' = 0$.  Its shifts $\Frob_{d,d'}^\oco$ are defined in the obvious way.  The phrase ``open Frobenius algebra'' (which seems to have originated in or around Dennis Sullivan, but I was not able to come up with an earliest reference) comes from the fact that the cohomology of an \emph{open} oriented manifold is a \emph{non-counital} shifted Frobenius algebra.

An immediate generalization of the Boardman--Vogt tensor product: if $P$ is an arrowed dioperad and $Q$ is a singly-colored dioperad, then $(P \boxtimes Q)(m,n;m',n') = P(m,n;m',n') \otimes Q(m+n;m'+n')$ is naturally arrowed.
Suppose now that $P$ and $Q$ are both arrowed.  There is a general notion of ``tensor product of colored dioperads'' that when fed two bicolored dioperads produces a four-colored dioperad; letting $\Arr$ denote the dioperad from section~\ref{section arrowed operad} so that arrowed dioperads are plain dioperads with a map from $\Arr$, we see that if $P$ and $Q$ are both arrowed dioperads, then $P \boxtimes Q$ receives a map from the dioperad $\Arr \boxtimes \Arr$ that parameterizes commuting squares.  

We will use instead two more interesting tensor products.  First, recall the tensor product of arrows from equation~(\ref{tensor product of arrows}):
$$ (A \overset f \to B) \quad \otimes \quad (X \overset p \to Y) \quad = \quad (A \otimes X) \overset {f \otimes p} \longto (B \otimes Y) $$
Suppose that $(A \overset f \to B)$ is an arrowed $P$-algebra and $(X \overset p \to Y)$ is an arrowed $Q$-algebra.  The arrowed dioperad acting on $(A \overset f \to B) \otimes (X \overset p \to Y)$ is 
$$ (P \boxtimes_\Arr Q)(m,n;m',n') = P(m,n;m',n') \otimes Q(m,n;m',n')$$
with arrowing $f\otimes p \in P(1,0;0,1) \otimes Q(1,0;0,1)$.  For example, if $P$ and $Q$ are plain dioperads, then $(P\boxtimes Q)^{\triangleleft,\str} \cong P^{\triangleleft,\str} \boxtimes_\Arr Q^{\triangleleft,\str}$.  

As discussed in section~\ref{section rotation}, the tensor product of arrows does not play well with rotation.  Our second tensor product of arrowed dioperads $P,Q$ is the plain dioperad $P \boxtimes_\Arr^\gamma Q = \gamma(P \boxtimes_\Arr Q)$.  Equation~(\ref{rotation and tensoring}) then implies:
\begin{equation} \label{rotation and tensoring of dioperads} P\boxtimes_\Arr^\gamma Q  \simeq \Theta P \boxtimes_\Arr^\gamma \Theta^{-1} Q.\end{equation}

\subsection{Example} \label{cone of an ideal}

  Let us investigate $\gamma(O^{\triangleleft,\str})$.  Combining equations~(\ref{strict ideal}) and~(\ref{cone of an arrowed dioperad}) gives
$$ \gamma O^{\triangleleft,\str}(k;k') = \bigoplus_{\substack{m,n'\in \{0,\dots,k\} \\ (m,n') \neq (0,0)}} \binom k m \binom{k'}{n'} O(k;k')[k'][-m-n'], $$
where here $\binom k m\binom{k'}{n'}$ means the $S_k \times S_{k'}$-module should be restricted to an $S_m \times S_{k-m} \times S_{k'-n'} \times S_{n'}$-module and then induced back.
The differential is nothing but a combination of the canonical maps connecting the $(m,n')$th summand to the $(m+1,n')$th and $(m,n'+1)$th summands.  If there were a $(0,0)$th summand, then the complex would be exact.  (This categorifies the fact that $\sum_{m=0}^k (-1)^m \binom k m = 0$).  In cohomology, therefore, $\gamma O^{\triangleleft,\str}(k;k')$ looks like a shifted copy of the would-be $(0,0)$th summand $O(k;k')[k']$, and there is a quasiisomorphism
$$ \gamma O^{\triangleleft,\str} \simeq O \boxtimes \Frob_{0,-1}.$$

\subsection{} \label{bar dual}

Let $P$ be a locally finite oco arrowed dioperad, and $\underline P$ the non-unital bicolored dioperad produced from $P$ by removing its identity elements and arrowing, as in section~\ref{section arrowed operad}.
The \define{bar dual} $\bD P$ of $P$ is the quasifree arrowed dioperad with generators $\underline P(m,n;m',n')^*[m+n'-2]$, where $V^*$ denotes the linear dual to the cochain complex $V$.  The differential on $\bD P$ is of course defined on generators, and has three terms:
\begin{enumerate} 
\item The linear differential on $\underline P^*$.
\item A term coming encoding the arrowing in $P$.  Recall that $\underline P$ has maps $f\circ : \underline P(m,n;m',n') \to \underline P(m,n;m'-1,n'+1)$ and $\circ f : \underline P(m,n;m',n') \to \underline P(m+1,n-1;m',n')$, and so the dual has maps $(f\circ)^*,(\circ f)^* : \underline P^* \to \underline P^*$.  (There are in fact many such maps depending on the choice of leaf at which to compose.)  Call the new arrowing of $\bD P$ ``$\phi$.''  Then the second term in the differential on $\bD P$ takes a generator $x\in \underline P(m,n;m',n')^*[m+n'-1]$ to $\sum \pm \phi \circ \bigl((f\circ)^*(x)\bigr) \pm \bigl((\circ f)^*x\bigr) \circ \phi$.  The sum is over leaves at which to do the composition and the signs depend on conventions that we leave implicit.
\item A term encoding the binary composition in $\underline P$.  There are various composition maps $\underline P(m_1,n_1;m_1',n_1') \otimes \underline P(m_2,n_2;m_2',n_2') \to \underline P(m_1 + m_2 - 1,n_1+n_2;m_1' + m_2'-1,n_1'+n_2')$ and $\underline P(m_1,n_1;m_1',n_1') \otimes \underline P(m_2,n_2;m_2',n_2') \to \underline P(m_1 + m_2,n_1+n_2-1;m_1' + m_2',n_1'+n_2-1')$, which in turn give maps $\underline P^* \to \underline P^* \otimes \underline P^*$.  We use these maps on generators.
\end{enumerate}
As usual, that the resulting derivation of $\bD P$ squares to zero is equivalent to the associativity of composition in $P$.  That the terms are all of the correct degree follows from simple combinatorics.

The second and third terms in the differential, encoding respectively the arrowing and the composition, can be combined by working not with $\underline P$ but with $\underline P \oplus \bK f$, which is again a non-unital bicolored dioperad.  Indeed, $\bD P$ is generated as a bicolored, rather than arrowed, dioperad by $(\underline P \oplus \bK f)^*[m+n'-2] = \underline P^*[m+n'-2] \oplus (\bK f)^*[1+1-2]$, and the differential on $\bD P$ can be seen as the linear differential plus a single term encoding binary composition in $\underline P \oplus \bK f$.  The arrowing of $\bD P$ is the dual basis vector to $f$ in $(\bK f)^*$.

The construction $\bD$ is called a \define{duality} because there is a quasiisomorphic surjection $\bD^2 P \overset\sim\epi P$.  This can be seen  readily with a bicolored version of the usual argument. It is worth mentioning that non-oco dioperads do not avoid the usual subtleties with bar duality. For example, $\bD$ extends natuarally to augmented but non-oco dioperads, but without some conditions it does not even preserve quasi-isomorphisms.

For an oco plain dioperad $O$, we let $\underline O$ denote the nonunital dioperad of operations in $O$ with total arity $\geq 3$ and $\bD O$ the quasifree dioperad generated by $\underline O^*[-1]$, with differential encoding binary composition in $\underline O$.  For example, $\bD (\Frob^\oco)$ is the dioperad of directed trivalent-and-higher trees with no sources or sinks; composition is by concatenation of trees and the differential takes a tree $\Gamma$ to the sum over all trees that can produce $\Gamma$ by contracting a single edge.  For any $O$, $\bD O$ is universal for there to be a map $\bD (\Frob^\oco) \to O \boxtimes \bD O$ --- the map sends the $(k;k')$th generator of $\bD(\Frob^\oco)$ to the canonical element in $O(k;k') \otimes O(k;k')^*[-1] \subseteq O \boxtimes \bD O$ ---
 in the sense that any map $\bD (\Frob^\oco) \to O \boxtimes Q$ factors through a map $\bD O \to Q$. 
 The corresponding well-known statement for operads is that if $O$ is an operad, $\bD O$ is universal for there to be a map $\bD \mathrm{Com}^\oco = L_\infty\langle 1\rangle \to O \boxtimes \bD O$.
  It follows in particular that 
 \begin{equation} \label{D and Frob} \bD(O \boxtimes \Frob_{d,d'}) \cong (\bD O) \boxtimes \Frob_{-d,-d'}.\end{equation}

\subsection{Theorem} \label{universal property of the dual}
\emph{Fix an oco arrowed dioperad $P$.  Then $\bD P$ is the universal arrowed dioperad equipped with a map $\bD \Frob_{0,1}^\oco \to P \boxtimes_\Arr^\gamma \bD P$.}

\subsection*{Proof}

We have $\Frob_{0,1}^\oco(k;k') = \bK[1-k']$, provided $kk' \geq 1$, and so $\bD \Frob_{0,1}^\oco$ is generated by $\underline{\Frob_{0,1}^\oco}^*(k;k')[-1] = \bK[k'-2]$ for $kk' \geq 2$.  Other than the degree shifts, $\bD \Frob_{0,1}^\oco$ is just the directed tree operad $\bD \Frob^\oco$: the differential takes a generator to the sum over two-vertex trees with the given arity.

Suppose that $Q$ is an arrowed dioperad with a dioperad map $\varphi: \bD \Frob_{0,1}^\oco \to P \boxtimes_\Arr^\gamma Q$.  Temporarily ignore all differentials.  Equation~(\ref{cone of an arrowed dioperad}) gives
$$ (P \boxtimes_\Arr^\gamma Q)(k;k') = \bigoplus_{\substack{m+n=k \\ m'+n'=k'}} \binom k m \binom{k'}{m'} P(m,n;m',n') \otimes Q(m,n;m',n') [m'-m]. $$
The map $\varphi$ consists of only the images of generators, and so for each fixed $(k;k')$ we have an $S_k \times S_{k'}$-equivariant map
$$ \bK[k'-2] \to \bigoplus_{\substack{m+n=k \\ m'+n'=k'}} \binom k m \binom{k'}{m'} P(m,n;m',n') \otimes Q(m,n;m',n') [m'-m]$$
or equivalently a degree-$0$ $S_k \times S_{k'}$-fixed element of
$$ \bigoplus_{\substack{m+n=k \\ m'+n'=k'}} \binom k m \binom{k'}{m'} P(m,n;m',n') \otimes Q(m,n;m',n') [2-m-n'],$$
where we have used $k'=m'+n'$.

A fixed point of a finite direct sum of modules is just a fixed point of each summand.  And fixed points of an induced module are fixed points of the module before induction.  Thus the map $\varphi$ is equivalent to, for each $(m,n;m',n')$, an $S_m \times S_n \times S_{m'} \times S_{n'}$-invariant element of $P(m,n;m',n') \otimes Q(m,n;m',n') [2-m-n']$.  This, in turn, is equivalent to an $S_m \times S_n \times S_{m'} \times S_{n'}$-equivariant map $P(m,n;m',n')^*[m+n'-2] \to Q(m,n;m',n')$.

The universal $Q$ is then freely generated by these $P(m,n;m',n')^*[m+n'-2]$s.  
  Restoring the differential to $\bD \Frob_{0,1}^\oco$ gives the stated differential for $\bD P$.
\qed\\[-4pt]

Combining Theorem~\ref{universal property of the dual} with equation~(\ref{rotation and tensoring of dioperads}) gives:

\subsection{Corollary}\label{rotation and duality}
\emph{For any oco arrowed dioperad $P$, $\bD \Theta P$ and $\Theta^{-1}\bD P$ are quasiisomorphic.} \qed

\subsection{} \label{quadratic}

A bicolored dioperad is \define{quadratic} if it is presented by generators and relations such that all relations are quadratic in the generators.  An arrowed dioperad is \define{quadratic} if it is quadratic as a bicolored dioperad, where the arrowing is included as a generator.  Given a collection $G$ of generators, let us write $\cF(G)$ for the free arrowed dioperad on $G$ (i.e.\ the free bicolored operad on $G \sqcup \{f\}$), and $\cF^{(w)}(G)$ the part of $\cF(G)$ of homogeneous ``polynomial'' degree $w$ in $G \sqcup \{f\}$. The quadratic relations are then a subcomplex $R \subseteq \cF^{(2)}(G)$, and a quadratic arrowed dioperad is presented as $P = \cF(G) / \langle R\rangle$, where $\langle R\rangle$ denotes the arrowed dioperad ideal generated by $R$.  A quadratic dioperad is \define{oco} if $G$ vanishes in arities $(m,n;m',n')$ with $(m+n)(m'+n') \leq 1$.

Suppose that $P = \cF(G) / \langle R\rangle$ is a locally finite oco quadratic arrowed dioperad.  The \define{quadratic dual} $P^!$ of $P$ is the quadratic arrowed dioperad with generators $G^!(m,n;m',n') = G(m,n;m',n')^*[m+n'-2]$ and relations $$R^\perp[m+n'-3] = \ker\Bigl(\cF^{(2)}(G^!) = \cF^{(2)}(G)^*[m+n'-3] \to R^*[m+n'-3]\Bigr).$$  

Quadratic duality is closely related to the bar duality of section~\ref{bar dual}.  There is an obvious surjection $\bD P \epi P^!$ sending the polynomial-degree-one piece of $P^*[m+n'-2]$ of $\bD P$ to the generators of $P^!$ and annihilating the rest of $P^*[m+n'-2]$; it is a map of arrowed dioperads because the derivatives of degree-two generators of $\bD P$ map to relations in $P^!$.
The quadratic dioperad $P$ is \define{Koszul} if the surjection $\bD P \epi P^!$ is a quasiisomorphism.  If $P$ is Koszul, so is $P^!$, since the quasiisomorphism $\bD^2 P \to P$ factors as $\bD^2 P \isom \bD P^! \to P$.

The notions of ``quadratic'' and ``Koszul'' make sense also for plain dioperads.  If $O$ is a plain quadratic dioperad with generators $G$ and relations $R$, $O^!$ is the quadratic dioperad with generators $G^*[-1]$ and relations $R^\perp[-2]$.

\subsection{Theorem} \label{duality commutes with ideals}  \emph{Suppose that $O$ is an oco Koszul quadratic dioperad with quadratic dual $O^!$.  Then the arrowed dioperad $O^{\triangleleft,\str}$ is quadratic and Koszul with dual $(O^!)^{\triangleleft,\str}$.}

\subsection*{Proof}
Let $O$ have generators $G$ and quadratic relations $R$.  
Then $O^{\triangleleft,\str}$ is generated by, for each generator $x\in G(k;k')$, a generator $x_A \in G^{\triangleleft,\str}(1,k-1;k',0)$ and a generator $x_B \in G^{\triangleleft,\str}(0,k;k'-1,1)$.  (Since $O$ is assumed oco, both $k$ and $k'$ are positive for all generators.  We continue to assume that $G(k;k')$ is closed for the $S_k \times S_{k'}$-action, in which case the complexes of $x_A$s and $x_B$s are, respectively, the restrictions of this action to $S_{k-1}\times S_{k'}$ and $S_k \times S_{k'-1}$.)
Other arities of generators can then be formed by composing $x_A$ and $x_B$ with the arrowing $f$.

The relations are as follows.  Each relation $r\in R$ leads to (generally four, but fewer in low arities) relations formed by coloring the edges of $r$ either ``A'' or ``B'' so that each of the two vertices in $r$ has either no A-inputs and exactly one B-output or exactly one A-input and no B-outputs.  There are also relations involving the arrowing $f$.  First, for each generator $x\in G(k;k')$, we declare the relation $f \circ x_A - x_B \circ f \in R^{\triangleleft,\str}(1,k-1;k'-1,1)$.  Second, suppose $x\in G(k;k')$ with $k\geq 2$.  Let $(12) \in S_k$ denote the first two leaves, and $y \circ (12)$ the result of applying this permutation to an element $y$ of an $S_k$-module.  The second relation then says that $(x\circ(12))_A \circ f - (x_A \circ f) \circ (12) \in R^{\triangleleft,\str}(2,k-2;k',0)$.  Similarly, if $k'\geq 2$, we declare $f \circ ((12)\circ x)_B - (12) \circ (f\circ x_B) \in R^{\triangleleft,\str}(0,k;k-2,2)$.
In terms of modules, what we mean is the following.  A priori, a term like ``$x_A \circ f$'' transforms in the module formed by inducing to $S_{k-2} \times S_2$ the restriction of the $S_k$-module $G(k;k')$ to $S_{k-2}$.  The relation ``$(x\circ(12))_A \circ f = (x_A \circ f) \circ (12)$'' declares that in fact it transforms in the module formed by directly restricting from an $S_k$-action to an $S_{k-2} \times S_2$-action.  In pictures (with composition from top to bottom and with $f$ denoted by a solid bullet), the relations in $O^{\triangleleft,\str}$ describing how the generators relate to the arrowing are:
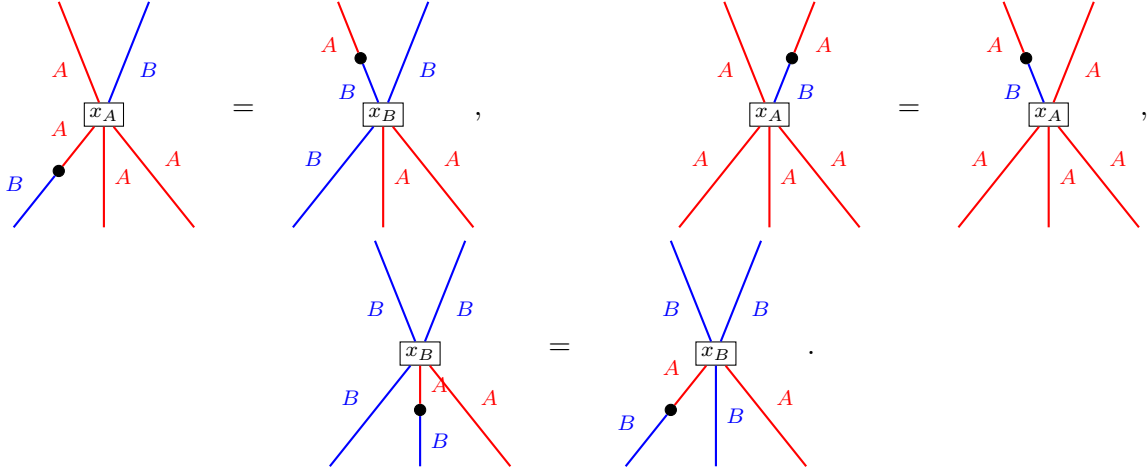
\begin{gather*}
\begin{tikzpicture}[scale=1.5,baseline=(x)]
  \path
    node[vertex] (x) {$\scriptstyle x_A$}
    ++(-.4,-.5) node[dot] (f) {}
  ;
  \draw [thick,red]
    (x) -- node[auto] {$\scriptstyle A$} +(-.4,1)
  ;
  \draw [thick,blue]
    (x) -- node[auto,swap] {$\scriptstyle B$} +(.4,1)
  ;
  \draw [thick,red]
    (x) -- node[auto,swap] {$\scriptstyle A$} (f) 
  ;
  \draw [thick,blue]
    (f) -- node[auto,swap] {$\scriptstyle B$} +(-.4,-.5)
  ;
  \draw  [thick,red]
    (x) -- node[auto] {$\scriptstyle A$} +(0,-1)
  ;
  \draw  [thick,red]
    (x) -- node[auto] {$\scriptstyle A$} +(.8,-1)
  ;
\end{tikzpicture}
  \quad = \quad
\begin{tikzpicture}[scale=1.5,baseline=(x)]
  \path
    node[vertex] (x) {$\scriptstyle x_B$}
    ++(-.2,.5) node[dot] (f) {}
  ;
  \draw [thick,blue]
    (x) -- node[auto,near end] {$\scriptstyle B$} (f)
  ;
  \draw [thick,red]
    (f) -- node[auto] {$\scriptstyle A$} +(-.2,.5)
  ;
  \draw [thick,blue]
    (x) -- node[auto,swap] {$\scriptstyle B$} +(.4,1)
  ;
  \draw [thick,blue]
    (x)  -- node[auto,swap] {$\scriptstyle B$} +(-.8,-1)
  ;
  \draw [thick,red]
    (x) -- node[auto] {$\scriptstyle A$} +(0,-1)
  ;
  \draw [thick,red]
    (x) -- node[auto] {$\scriptstyle A$} +(.8,-1)
  ;
\end{tikzpicture}
,\hspace{1in}
\begin{tikzpicture}[scale=1.5,baseline=(x)]
  \path
    node[vertex] (x) {$\scriptstyle x_A$}
    ++(.2,.5) node[dot] (f) {}
  ;
  \draw [thick,red]
    (x) -- node[auto] {$\scriptstyle A$} +(-.4,1)
  ;
  \draw [thick,blue]
    (x) -- node[auto,swap,near end] {$\scriptstyle B$} (f) 
  ;
  \draw [thick,red]
    (f) -- node[auto,swap] {$\scriptstyle A$} +(.2,.5)
  ;
  \draw[thick,red]
    (x) -- node[auto,swap] {$\scriptstyle A$} +(-.8,-1)
  ;
  \draw [thick,red]
    (x) -- node[auto] {$\scriptstyle A$} +(0,-1)
  ;
  \draw [thick,red]
    (x) -- node[auto] {$\scriptstyle A$} +(.8,-1)
  ;
\end{tikzpicture}
  \quad = \quad
\begin{tikzpicture}[scale=1.5,baseline=(x)]
  \path
    node[vertex] (x) {$\scriptstyle x_A$}
    ++(-.2,.5) node[dot] (f) {}
  ;
  \draw [thick,red]
    (x) -- node[auto,swap] {$\scriptstyle A$} +(.4,1)
  ;
  \draw [thick,blue]
    (x) -- node[auto,near end] {$\scriptstyle B$} (f) 
  ;
  \draw [thick,red]
    (f) -- node[auto] {$\scriptstyle A$} +(-.2,.5)
  ;
  \draw[thick,red]
    (x) -- node[auto,swap] {$\scriptstyle A$} +(-.8,-1)
  ;
  \draw [thick,red]
    (x) -- node[auto] {$\scriptstyle A$} +(0,-1)
  ;
  \draw [thick,red]
    (x) -- node[auto] {$\scriptstyle A$} +(.8,-1)
  ;
\end{tikzpicture}
,\\
\begin{tikzpicture}[scale=1.5,baseline=(x)]
  \path
    node[vertex] (x) {$\scriptstyle x_B$}
    ++(0,-.5) node[dot] (f) {}
  ;
  \draw [thick,blue]
    (x) -- node[auto] {$\scriptstyle B$} +(-.4,1)
  ;
  \draw [thick,blue]
    (x) -- node[auto,swap] {$\scriptstyle B$} +(.4,1)
  ;
  \draw [thick,red]
    (x) -- node[auto] {$\scriptstyle A$} (f) 
  ;
  \draw[thick,blue]
    (f) -- node[auto] {$\scriptstyle B$} +(0,-.5)
  ;
  \draw [thick,blue]
    (x) -- node[auto,swap] {$\scriptstyle B$} +(-.8,-1)
  ;
  \draw [thick,red]
    (x) -- node[auto] {$\scriptstyle A$} +(.8,-1)
  ;
\end{tikzpicture}
  \quad = \quad
\begin{tikzpicture}[scale=1.5,baseline=(x)]
  \path
    node[vertex] (x) {$\scriptstyle x_B$}
    ++(-.4,-.5) node[dot] (f) {}
  ;
  \draw [thick,blue]
    (x) -- node[auto] {$\scriptstyle B$} +(-.4,1)
  ;
  \draw [thick,blue]
    (x) -- node[auto,swap] {$\scriptstyle B$} +(.4,1)
  ;
  \draw[thick,red]
    (x) -- node[auto,swap] {$\scriptstyle A$} (f) 
  ;
  \draw[thick,blue]
    (f) -- node[auto,swap] {$\scriptstyle B$} +(-.4,-.5)
  ;
  \draw [thick,blue]
    (x) -- node[auto] {$\scriptstyle B$} +(0,-1)
  ;
  \draw [thick,red]
    (x) -- node[auto] {$\scriptstyle A$} +(.8,-1)
  ;
\end{tikzpicture}.
\end{gather*}

Inspection then reveals that $(O^{\triangleleft,\str})^!$ and $(O^!)^{\triangleleft,\str}$ are canonically isomorphic.  
For any oco dioperad $O$, the definitions directly give $\bD(O^{\triangleleft,\str}) \cong (\bD O)^\triangleleft$ and  Theorem~\ref{O ideal presents O strict ideal} says that the canonical surjection $ (\bD O)^\triangleleft \epi (\bD O)^{\triangleleft,\str} $ is a quasiisomorphism.  Koszulity of $O$ gives a surjective quasiisomorphism $\bD O \overset\sim\epi O^!$ and hence $ (\bD O)^{\triangleleft,\str} \overset\sim\epi (O^!)^{\triangleleft,\str}.$
All together, we see that the surjection
$$ \bD(O^{\triangleleft,\str}) \cong (\bD O)^\triangleleft \overset\sim\epi (\bD O)^{\triangleleft,\str} \overset\sim\epi (O^!)^{\triangleleft,\str} \cong (O^{\triangleleft,\str})^!$$
is a quasiisomorphism. \qed\\[-4pt]

Combining this with Theorem~\ref{O ideal justification} and Corollary~\ref{rotation and duality} gives:

\subsection{Corollary} \emph{Let $O$ be an oco Koszul operad.  The $\infty$ version of \define{extension} of $O$-algebras, completing the triangle whose other two sides are ``$\infty$-morphism'' and ``$\infty$-ideal'' of $\infty$-$O$-algebras, can be presented by the quasifree arrowed operad $\bD\bigl((O^!)^{\to,\str}\bigr)$.} \qed

\subsection{} \label{derived geometry}

Let $\Com$ denote the operad parameterizing unital commutative algebras and $\Com^\oco$ the operad parameterizing nonunital commutative algebras, so that $\Com(k;k') = \bK$ if $k'=1$ and $0$ otherwise, and $\Com^\oco(k;k') = \bK$ if $k'=1$ and $k\geq 1$ and $0$ otherwise.  $\Com^\oco$ is known to be Koszul; its quadratic dual is $(\Com^\oco)^! = \mathrm{Lie}\langle 1\rangle$.

A \define{strict dg affine scheme} is nothing but a $\Com$-algebra; we will write $\Spec(A)$ when we are thinking of a $\Com$-algebra $A$ as a dg affine scheme.  Strict morphisms of dg affine schemes are opposite to strict morphisms of $\Com$-algebras.  The corresponding notions of \define{homotopy dg affine scheme} and \define{$\infty$-morphism} thereof are achieved by replacing $\Com$ by some quasifree resolution $\mathrm{hCom}$ and using the arrowed operad $\mathrm{hCom}^\to$.  A \define{pointed strict dg affine scheme} is a map $\Spec(\bK) \to \Spec(A)$ of dg affine schemes, or equivalently a map $A \to \bK$ of $\Com$-algebras.  Such data is equivalent to giving $\ker(A \to \bK)$ the structure of a $\Com^\oco$-algebra, and so a minimal resolution of the notion of ``$\infty$-morphism of pointed dg affine schemes'' is afforded by the arrowed operad $\bD(\mathrm{Lie}\langle 1\rangle)^\to$.

A \define{pointed infinitesimal dg manifold} is a $\bD(\Com^\oco) = L_\infty\langle 1\rangle$-algebra.  Such data consists of: a cochain complex $(X,\partial_X)$ together with, for each $k\geq 2$, a map $\mu_k: \operatorname{Sym}^k(X) \to X[1]$.  Call $\partial_X = \mu_1 $ and consider the sum $\mu = \sum_{k\geq 1} \mu_k$.  The idea is to consider the graded vector space $X$ as a linear chart (centered at the pointing) for the corresponding infinitesimal manifold and the sum $\mu$ as a vector field on $X$; the differential in $\bD(\Com^\oco)$ says precisely that $\mu$ is cohomological.  (If $X$ is finite-dimensionsonal, then the algebra of functions on $X$ is the completed symmetric algebra $\prod_k \operatorname{Sym}^k(X^*)$ and the vector field is $\sum_k \mu_k^*(x) \frac\partial{\partial x} : X^* \to \operatorname{Sym}^k(X^*)$, at least up to some convention-dependent $k!$s.  A ``manifold'' is anywhere that you can do differential calculus; power series algebras over a field of characteristic $0$ certainly suffice.)  A morphism of pointed infinitesimal dg manifolds is an $\infty$-morphism of $\bD(\Com^\oco)$-algebras.

Suppose that $\Spec(A)$ is a strict dg affine scheme and $X$ is a pointed infinitesimal dg manifold.  Then the \define{mapping space} $\maps(\Spec(A),X)$ is the pointed infinitesimal dg manifold $A \otimes X$, which is given a $\bD(\Com^\oco)$-algebra structure using the canonical isomorphism $\Com \boxtimes \bD(\Com^\oco) \cong \bD(\Com^\oco)$.  This is reasonable because the underlying graded vector space of $X$ is supposed to be a linear chart for the pointed infinitesimal dg manifold and because we should have $A = \maps(\Spec(A),\bK)$.  If $\Spec(A)$ is merely a homotopy dg affine scheme, abstract nonsense of homotopical algebra assures that $\maps(\Spec(A),X) = A\otimes X$ still carries a $\bD(\Com^\oco)$-algebra structure, canonical up to a contractible space of choices depending on the chosen model $\mathrm{hCom}$ of ``homotopy commutative.''  Indeed, $\bD(\Com^\oco)$, being quasifree, is cofibrant for the model structure on operads in which surjections are fibrations, and there is a quasiisomoprhic surjection $\mathrm{hCom} \otimes \bD(\Com^\oco) \overset\sim\epi \bD(\Com^\oco)$.

A main idea of \cite{PoissonAKSZ} was to extend this well-known story from operads to dioperads.  A dg affine scheme $\Spec(A)$ is \define{$d$-oriented} if the $\Com$-algebra structure on $A$ is extended to a $\Frob_{0,d}$-algebra structure.  The idea is to think of the counit $i \in \Frob_{0,d}(1;0) = \bK[d]$ as giving an \define{integration} map $\int: A \to \bK[-d]$.  A \define{$\Pois_d$-structure} on a pointed infinitesimal manifold $X$ is an extension of the $\bD(\Com^\oco)$-structure to a $\bD(\Frob_{0,d}^\oco)$-structure.  The idea is that the generators in $\bD(\Frob_{0,d}^\oco)$ with arity $(k;2)$ together form the Poisson bivector field $\pi_2$ on $X$ just as the generators of arity $(k;1)$ formed a cohomological vector field $\mu = \pi_1$.  The Poisson bivector field $\pi_2$ is not itself a strict Poisson bivector field, but the trivector field $\pi_3$ formed from the generators of arity $(k;3)$ provides a ``Jacobiator'' for  $\pi_2$, and in general the polyvector fields $\pi_{k'}$ satisfy the rules of an $L_\infty\langle 1-d\rangle$-algebra.  The numbering is such that $\Pois_1$ corresponds to usual Poisson (``Pois-un'') with the bivector field $\pi_2$ in degree $0$.

We observed in section~\ref{bar dual} that $\bD(O \boxtimes \Frob_{0,d}) \cong (\bD O) \boxtimes \Frob_{0,-d}$.  It follows that:

\subsection{Theorem \cite{PoissonAKSZ}} \label{poisson aksz construction}
\emph{If $\Spec(A)$ is a $d$-oriented dg affine scheme and $X$ is a $\Pois_{d'}$ pointed infinitesimal manifold, then the pointed infinitesimal manifold $\maps(\Spec(A),X)$ is naturally $\Pois_{d'-d}$.} \qed\\[-4pt]

In fact, it suffices for $A$ to be merely $\Frob_{0,d}^\oco$.  The unit and counit are required to extend to the case when $X$ is, respectively, not pointed or ``curved.''  These versions require a curved variation of Koszul duality (see e.g.\ \cite{MR2993002}) and will not be described in this paper.

\subsection{} \label{d orientation and coiso defns}

Theorem~\ref{poisson aksz construction} is a Poisson generalization of the  AKSZ construction due to \cite{MR1432574}; see also \cite{MR3257656,MR3090262}.  Consideration of the symplectic case suggests that there should be a ``relative'' version of Theorem~\ref{poisson aksz construction}.  Specifically, suppose that $\Spec(C) \to \Spec(B)$ is a map of dg affine schemes where $\Spec(C)$ is $(d-1)$-oriented and $\Spec(B)$ is ``$d$-oriented relative to $\Spec(C)$'' --- we will give a precise definition momentarily, but the motivation comes from the case when $B = \H^\bullet(M)$ for $M$ an oriented $d$-dimensional manifold with boundary and $C = \H^\bullet(\partial M)$ --- and suppose that $Y \mono Z$ is a ``coisotropic submanifold'' of a $\Pois_d$-manifold $Z$.  Then the interesting thing to  study are those maps $\Spec(B) \to Z$ that when restricted to $\Spec(C)$ land in $Y$.

We now propose precise definitions.  Suppose that $M$ is an oriented manifold and $\partial M$ its boundary.  Then there is an exact triangle 
$\H^\bullet(\partial M)[-1] \to \H^\bullet(M;\partial M) \to \H^\bullet(M) \to \H^\bullet(\partial M)$
where $\H^\bullet(M;\partial M)$ are the relative cohomology groups of $M$ relative to $\partial M$.
The restriction map $\H^\bullet(M) \to \H^\bullet(\partial M)$ is a homomorphism of $\Com$-algebras.  Because $\partial M$ is closed and oriented, $\H^\bullet(\partial M) \cong \H_\bullet(\partial M)[1-d]$ and so is a coalgebra; relative Poincar\'e duality says that $\H^\bullet(M;\partial M)[1] \cong \H_\bullet(M)[1-d]$ is also a coalgebra; the ``corestriction'' map $\H_\bullet(\partial M) \to \H_\bullet(M)$ is a coalgebra homomorphism.  Finally,  the arrow $\H^\bullet(M;\partial M) \overset i \to \H^\bullet(M)$ is easily checked to be a $(\Frob_{0,d})^{\triangleleft,\str}$-algebra.  We therefore declare:

\subsection*{Definition} A map $\Spec(C) \to \Spec(B)$ of dg affine schemes is \define{relative $d$-oriented} if the $\Com^\to$-algebra structure on the arrow $B \to C$ is extended to $\Theta\bigl( \hFrob_{0,d}^\triangleleft\bigr)$, where $\hFrob_{0,d}$ is any quasifree resolution of $\Frob_{0,d}$.\\[-4pt]

Theorem~\ref{O ideal justification} assures that $\Com^\to \simeq \Theta\Com^\triangleleft$, so this definition is sensical.  
Note that a relative orientation on $\Spec(B) \to \Spec(A)$ includes the data of a $(d-1)$-orientation on $\Spec(B)$ by Example~\ref{cone of an ideal}, which implies that $\beta\Theta\hFrob_{0,d}^\triangleleft \simeq \Frob_{0,d-1}$.

Dually, what should be a ``coisotropic sub'' $Y$ of a Poisson pointed infinitesimal dg manifold~$Z$?  
Certainly we should have a map $Y \to Z$ of pointed infinitesimal dg manifolds.  Such data is parameterized by $\bD(\Com^\oco)^\to$.
  Suppose that $Y$ and $Z$ were not pointed infinitesimal dg manifolds but ordinary manifolds.  Then to say that $Z$ is Poisson would be to say that its algebra of functions $\cO(Z)$ is Lie (in a way compatible with the commutative structure), and to say that $Y \mono Z$ is coisotropic would be to say that $K = \ker\bigl(\cO(Z) \to \cO(Y)\bigr)$ should be a Lie subalgebra of $\cO(Z)$.
The Lie structure on $\cO(Z)$ corresponds in the pointed infinitesimal dg case to a $\mathrm{Lie}^\op$ structure on $Z$, or more accurately a $\bD(\Com^{\oco,\op}\langle- d\rangle)$-structure.  Thus the homomorphism $K \to \cO(Z)$ of Lie algebras should correspond to a homomorphism $Z \to K'$ of $\bD(\Com^{\oco,\op}\langle -d\rangle)$-algebras.  

Setting $X = K'[-1]$, we find ourselves with an exact triangle $Z[-1] \to X \to Y \to Z$ where $Z[-1] \to X$ is a $\bD(\Com^{\oco,\op}\langle 1-d\rangle)^\to$-algebra and $Y\to Z$ is a $\bD(\Com^\oco)^\to$-algebra.  The arrow $X \to Y$ is both a $\bD(\Com^\oco)$-ideal and a $\bD(\Com^{\oco,\op}\langle 1-d\rangle)$-coideal, and so it is reasonable to demand that it should be a $\bD(\Frob_{0,d-1}^\oco)^\triangleleft$-algebra.

\subsection*{Definition}
A map $Y\to Z$ of pointed infinitesimal dg manifolds is a \define{coisotropic in $\Pois_d$} if the $\bD(\Com^\oco)^\to$-algebra structure on the arrow $Y \to Z$ is extended to $\Theta \bigl(\bD(\Frob_{0,d-1}^\oco)^\triangleleft\bigr)$.
\\[-4pt]

Other definitions of ``coisotropic'' are discussed in \cite{Safronov2015,MelaniSafranov2016,Safronov2016}. 
As discussed in those papers, there are already many definitions of ``coisotropic'' in derived geometry; only very recently have many of them been shown to be equivalent.
Checking that the definition from this paper matches the others will be the subject of future work.

We are now equipped to prove the following ``relative''  version of Theorem~\ref{poisson aksz construction}, which constitutes the algebraic half of the ``Poisson AKSZ construction with coisotropic boundary conditions'':

\subsection{Theorem}\label{AKSZ construction with boundary}

\emph{Let $\Spec(C)$ be a $(d-1)$-oriented affine dg scheme and $\Spec(C) \to \Spec(B)$ a map of affine dg schemes equipped with a relative $d$-orientation.  Let $Z$ be a $\Pois_{d'}$ pointed dg infinitesimal manifold and $Y \to Z$ a map of pointed dg infinitesimal manifolds equipped with a coisotropic structure.  Note that Theorem~\ref{poisson aksz construction} gives the pointed infinitesimal dg manifold $\maps(\Spec(C),Z)$ a $\Pois_{d'-d+1}$-structure.  The pointed infinitesimal dg manifolds $\maps(\Spec(B),Z)$, and $\maps(\Spec(C),Y)$ are coisotropics in $\maps(\Spec(C),Z)$. The mapping space
$$
\maps(\Spec(B),Z) \underset{\maps(\Spec(C),Z)}{\times^h} \maps(\Spec(C),Y)
$$
of maps $\Spec(B) \to Z$ whose restriction lands to $\Spec(C)$ lands in $Y$
is $\Pois_{d'-d}$.  (The ``$h$'' emphasizes that this is a homotopy fibered product.)  The space $\maps(\Spec(B),Y)$ is coisotropic therein.}
%\\[-4pt]

\subsection*{Proof}

That the arrow $\bigl(\maps(\Spec(C),Y) \to \maps(\Spec(C),Z)\bigr) = \bigl(C \otimes Y \to C \otimes Z\bigr)$ is coisotropic follows from the fact that the tensor product of the arrowed dioperad $\Theta\bigl(\bD(\Frob_{0,d'-1}^\oco)^\triangleleft\bigr)$ presenting the notion of ``coisotropic'' with the plain dioperad $\Frob_{0,d-1}$ describing the orientation on $\Spec(B)$ is
$$ \Frob_{0,d-1} \boxtimes \Theta \bigl(\bD(\Frob_{0,d'-1}^\oco)^\triangleleft\bigr) = \Theta \bigl(\bD(\Frob_{0,d'-d}^\oco)^\triangleleft\bigr).$$
The statement that $\maps(\Spec(B),Z) \to \maps(\Spec(C),Z)$ is coisotropic follows from the Koszul dual argument.

Extend the arrows $B \to C$ and $Y \to Z$ to exact triangles $C[-1] \to A \to B \to C$ and $Z[-1] \to X \to Y \to Z$.
Corollary~\ref{rotation and duality}  implies that the notions ``coisotropic sub'' and ``relative orientation'' are almost Koszul dual.  
After rotating to  line up the Koszul duality, we are led to consider the following \define{relative AKSZ exact Star of David}:

\begin{equation} \label{AKSZ Star}
\begin{tikzpicture}[baseline=(middle),thick]
  \path 
  (150:4.5) node (A) {$A$}  (30:4.5) node (B) {$B$} (-90:4.5) node (C1) {$C[-1]$} ++(0,18pt) node (C) {$C$} 
  (-150:4.5) node (X) {$X$}  (-30:4.5) node (Y) {$Y$} (90:4.5) node (Z1) {$Z[-1]$}  ++(0,-18pt) node (Z) {$Z$}
  (0,0) coordinate (middle);
  \draw[->] (Z1) -- node[fill=white,inner sep=1pt,pos=.46] {$\scriptstyle (\bD\Com^{\oco,\op}\langle d'\rangle)^\to\!$} (X);
  \draw[->] (Y) -- node[fill=white,inner sep=1pt,pos=.6] {$\scriptstyle (\bD\Com^\oco)^\to$} (Z);
  \draw[->] (X) -- node[fill=white,inner sep=1pt,pos=.48] {$\scriptstyle (\bD\mathrm{Frob}_{0,d'-1}^\oco)^\triangleleft$} (Y);
  \draw[->] (C1) -- node[fill=white,inner sep=2pt,pos=.48] {$\scriptstyle \mathrm{hCom}^\op\langle d-1\rangle^\to$} (A);
  \draw[->] (B) -- node[fill=white,inner sep=2pt,pos=.57] {$\scriptstyle \mathrm{hCom}^\to$} (C);
  \draw[->] (A) -- node[fill=white,inner sep=1pt] {$\scriptstyle \mathrm{hFrob}_{0,d}^\triangleleft$} (B);
\end{tikzpicture}
\end{equation}

Focus on the arrow $A \to B$, which up to homotopical replacements is a $\Frob_{0,d}^\triangleleft$-algebra, and the arrow $X \to Y$, which is a $\bD(\Frob^\oco_{0,d'-1})^\triangleleft$-algebra.  Equation~(\ref{strict ideal}) immediately implies that, for any plain oco dioperad $O$, $\Frob_{0,d}^{\triangleleft,\str} \boxtimes_\Arr O^{\triangleleft,\str} \cong (\Frob_{0,d} \boxtimes O)^{\triangleleft,\str}$, and along with equation~(\ref{D and Frob}) we conclude that the tensor product of arrows $(A\otimes X) \to (B \otimes Y)$ is a $\bD(\Frob_{0,d'-d-1})^\triangleleft$-algebra.
It follows from Example~\ref{cone of an ideal} that $\cone\bigl((A\otimes X) \to (B \otimes Y)\bigr)$ is a $\bD\Frob_{0,d'-d}$-algebra,  i.e.\ a $\Pois_{d'-d}$ pointed infinitesimal dg manifold, and that $B \otimes Y$ is a coisotropic therein.

But $\cone\bigl((A\otimes X) \to (B \otimes Y)\bigr)$ is precisely the canonical ``tensor product'' cochain complex of the star~(\ref{AKSZ Star}) described in section~\ref{section rotation}.  It follows that
$$ \cone\bigl((A\otimes X) \to (B \otimes Y)\bigr) \simeq \cone\bigl( (C[-1] \otimes Y) \to (A \otimes Z) \bigr) $$
and unpacking $A \simeq \cone(B \to C)[-1]$ shows that this is quasiisomorphic to
\\[6pt] \mbox{} \hfill
$ \simeq (B \otimes Z) \underset{(C\otimes Z)}{\times^h} (C \otimes Y) = \maps(\Spec(B),Z) \underset{\maps(\Spec(C),Z)}{\times^h}\maps(\Spec(C),Y).$ \qed

\subsection{} \label{quasilocal operators}
For the remainder of this paper we work over the ground field $\bK = \bR$.
One of the main results of \cite{PoissonAKSZ} was the construction, for any closed oriented manifold $M$ of dimension $d$, of a canonical $\Frob_{0,d}^\oco$-algebra structure on the de Rham complex $\Omega^\bullet(M)$ satisfying a locality-type condition called \define{quasilocality} that is important for (classical and quantum) field theory.  Quasilocality assures that the $\Pois_0$ structure on $\maps(M_{\mathrm{dR}},X) = \maps(\Spec(\Omega^\bullet(M)),X)$, where $X$ is a $\Pois_d$ infinitesimal manifold, satisfies a weak sheaflike condition: as the ``energy scale'' of the theory increases, the ``support'' of the bracket of ``observables'' can be made as close to the support of the observables as desired.  Similar ideas are vital in \cite{CG}.
We end this paper by performing a similar construction for manifolds with boundary.  We use the improved version of ``quasilocality'' from~\cite{me-Qloc}.  

Let $M$ and $N$ be compact manifolds, possibly with boundary.  For convenience we assume them both oriented --- otherwise the following discussion must be decorated with twists by orientation bundles which merely gum up the notation.  Given a submanifold $Y \mono M$, $\Omega^\bullet(M;Y)$ denotes the relative de Rham complex  of forms on $M$ that restrict to $0$ on $Y$, and $\H^\bullet(M;Y)$ the corresponding relative cohomology groups.  The orientation picks out an isomorphism $\H^\bullet(M)^* \cong \H^\bullet(M;\partial M)[\dim M]$, and so
\begin{gather*}
  \H^\bullet \hom\bigl(\Omega^\bullet(M),\Omega^\bullet(N)\bigr) \cong \H^\bullet\bigl(M \times N; \partial M \times N\bigr) [\dim M] \\ 
  \H^\bullet \hom\bigl(\Omega^\bullet(M;\partial M),\Omega^\bullet(N)\bigr) \cong \H^\bullet\bigl(M \times N\bigr) [\dim M] \\
  \H^\bullet \hom\bigl(\Omega^\bullet(M),\Omega^\bullet(N;\partial N)\bigr) \cong \H^\bullet\bigl(M \times N;\partial M \times N\bigr) [\dim M] \\
  \H^\bullet \hom\bigl(\Omega^\bullet(M;\partial M),\Omega^\bullet(N;\partial N)\bigr) \cong \H^\bullet\bigl(M \times N; M \times \partial N\bigr) [\dim M] 
\end{gather*}
where in all cases $\hom$ means the complex of continuous linear maps.  By identifying a linear map with its integral kernel, we can in fact think of these complexes $\hom(\Omega\dots)$ as complexes of (singular) de Rham forms on $M \times N$ with the given boundary conditions.  An element $\psi$ of any of these complexes \define{avoids} a point $(x,y) \in M \times N$ if there are small neighborhoods $U \ni x$ and $V \ni Y$ such that for any $\omega$ supported entirely in $U$, $\psi(\omega)$ vanishes in $V$.  More generally, $\psi$ \define{avoids} a submanifold $Y \subseteq M \times N$ if it avoids every point in $Y$.  The \define{support} $\supp(\psi) \subseteq M \times N$  is the set of all points that $\psi$ does not avoid.

Let $\epsilon$ be a parameter ranging in $\bR_{>0}$.  Physically, $\epsilon$ should be thought of as a ``length scale,'' so that $\epsilon^{-1}$ is an ``energy scale.''
A \define{homotopy-constant family} in any of the complexes $\hom(\Omega\dots)$ above is an expression of the form $\psi(\epsilon) + \phi(\epsilon) \d \epsilon$, where $\psi,\phi$ depend smoothly on $\epsilon$.  The differential on homotopy-constant families is $\partial(\psi(\epsilon)) = [\d,\psi(\epsilon)] + \bigl(\frac{\partial \psi}{\partial \epsilon} + [\d,\phi(\epsilon)]\bigr) \d\epsilon$, where $[\d,-]$ denotes the differential on $\hom(\Omega\dots)$.
  The name ``homotopy-constant''comes from the fact that the inclusion of actually-constant among homotopy-constant families is a quasiisomorphism.  Indeed, $\psi(\epsilon) + \phi(\epsilon) \d \epsilon$ is closed if and only if for every $\epsilon_1,\epsilon_2 \in \bR_{>0}$, $\psi(\epsilon_1) - \psi(\epsilon_2) = \bigl[\d,\int_{\epsilon_1}^{\epsilon_2}\phi(\epsilon) \d \epsilon\bigr]$, so that the values of $\psi$ are homotopic by a prescribed homotopy.
  For the appropriate completed tensor product, the complex of homotopy-constant families is $\hom(\Omega\dots) \otimes \Omega^\bullet(\bR_{>0})$, and we will use this formula rather than introducing a new term.  
  Because $\Omega^\bullet(\bR_{>0})$ is a strict $\Com$-algebra, homotopy-constant families of operations compose without trouble: one composes in the ``operations'' direction and multiplies in the ``$\epsilon$'' direction.

Let $L \mono M \times N$ be a compact oriented submanifold.  A homotopy-constant family $\psi(\epsilon) + \phi(\epsilon) \d \epsilon$ is \define{near $L$} if for any open set $U \supseteq L$, there is a cut-off $\epsilon_U \in \bR_{>0}$ such that for all $\epsilon < \epsilon_U$, the linear maps $\psi(\epsilon)$ and $\phi(\epsilon)$ are supported in $U$.  (This should properly be called ``ultravioletly near $L$,'' since one could also consider ``infrared'' behavior of $\psi$ and $\phi$ as $\epsilon \to \infty$.)
As explained in~\cite{me-Qloc}, requiring a homotopy-constant family to be near $Y$ is equivalent to placing a support condition on the integral kernel of $\psi(\epsilon) + \phi(\epsilon) \d \epsilon$ on $\bR_{>0} \times M \times N$:
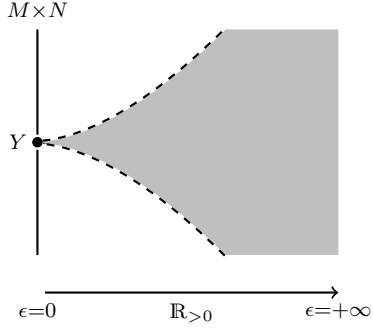
$$
  \begin{tikzpicture}[baseline=(diag)]
    \path (0,0) node[circle,draw,fill,inner sep=1.25pt] (diag) {} node[anchor=east] {$\scriptstyle Y$} ;
    \draw[thick] (0,.1) -- (0,1.5) node[anchor=south] {$\scriptstyle M\times N$} ;
    \draw[thick] (0,-.1) -- (0,-1.5);
    \draw[thick,->] (0,-2) node[anchor=north] {$\scriptstyle \epsilon = 0$} (.1,-2) -- node[auto,swap] {$\scriptstyle \bR_{>0}$} (4,-2) node[anchor=north] {$\scriptstyle \epsilon = +\infty$};
    \draw[ultra thick,dashed] (diag) .. controls +(1,.1) and +(-.5,-.5) .. (2.5,1.5);
    \draw[ultra thick,dashed] (2.5,-1.5) .. controls  +(-.5,.5) and +(1,-.1) .. (diag);
    \path[fill=black!25] (diag) .. controls +(1,.1) and +(-.5,-.5) .. (2.5,1.5) -- (4,1.5) -- (4,-1.5) -- (2.5,-1.5) .. controls  +(-.5,.5) and +(1,-.1) .. (diag);
  \end{tikzpicture}
$$
Moreover, a Thom-type isomorphism identifies the cohomology of the space of near-$L$ families with the relative cohomology $\H^\bullet(L;Z)[\dim L - \dim N]$, where $Z$ depends on which of the $\hom(\Omega\dots)$ complexes is being considered.

Consider the arrowed dioperads $\End\bigl( \Omega^\bullet(M;\partial M)\hookrightarrow \Omega^\bullet(M)\bigr)$ of \define{bulk-bulk operations on $M$} and $E(M) = \End\bigl( \Omega^\bullet(M;\partial M)\hookrightarrow \Omega^\bullet(M)\bigr)\otimes \Omega^\bullet(\bR_{>0})$ of homotopy-constant families of bulk-bulk operations.  There is a quasiisomorphism $\End \simeq E$, and so in cohomology $\H^\bullet(E(M)) \cong \End\bigl( \H^\bullet(M;\partial M) \to \H^\bullet(M)\bigr)$.
Let $\diag : M \to M^k$ denote the diagonal embedding.
Inside of $E(M)$ we define a subdioperad $\Qloc(M)$ of \define{quasilocal bulk-bulk operations} by
$$ \Qloc^\triangleleft(M) (m,n;m',n') = \bigl\{ \psi \in E(M)(m,n;m',n') : \text{$\psi$ is near $\diag(M) \mono M^{m+n+m'+n'}$}\bigr\} $$
That $\Qloc^\triangleleft(M)$ is closed under dioperadic composition follows from the triangle inequality.

Similarly, inside the arrowed dioperad $\End\bigl(\Omega^\bullet(M) \to \Omega^\bullet(\partial M)\bigr) \otimes \Omega^\bullet(\bR_{>0})$ of homotopy-constant families of \define{bulk-boundary operations}, we can define a subdioperad
$$\Qloc^\to(M)(m,n;m',n') = \begin{cases} \psi\text{ is near }\diag(\partial M) \mono M^{m+m'} \times (\partial M)^{n+n'} , & n+n' > 0,\\
 \psi\text{ is near }\diag( M) \mono M^{m +m'} , & n=n' = 0 \end{cases} $$
of \define{quasilocal bulk-boundary operations}.  For any neighborhood $U \supseteq \partial M$ in $M$, one can find a deformation retraction of $\cone\bigl( \Omega^\bullet(M) \to \Omega^\bullet(\partial M)\bigr)[1]$ onto $\Omega^\bullet(M;\partial M)$ that is supported in $U$.  By shrinking $U$ with $\epsilon$, one can build a quasiisomorphism $\Theta \Qloc^\to(M) \simeq \Qloc^\triangleleft(M)$.

\subsection{Theorem} \label{quasilocal frobenius structure}

\emph{Let $\mathrm{hFrob}_{0,d}^\oco$ denote any quasifree resolution of $\Frob_{0,d}^\oco$.  Let $M$ be a compact oriented $d$-dimensional manifold with boundary.
Up to a contractible space of choices, there is a unique map of arrowed dioperads $(\mathrm{hFrob}_{0,d}^\oco)^\triangleleft \to \Qloc^\triangleleft(M)$ such that the induced map on cohomology $(\Frob_{0,d}^\oco)^{\triangleleft,\str} = \H^\bullet(\mathrm{hFrob}_{0,d}^\oco)^\triangleleft \to \H^\bullet\Qloc(M^\triangleleft) \to \H^\bullet(E(M)) = \End\bigl(\H^\bullet(M;\partial M) \to \H^\bullet(M)\bigr)$ is the canonical $(\Frob_{0,d}^\oco)^{\triangleleft,\str}$-algebra structure on the arrow $\H^\bullet(M;\partial M) \to \H^\bullet(M)$ determined by Poincar\'e duality.}\\[-4pt]

It follows that the arrow $\Omega^\bullet(M) \to \Omega^\bullet(\partial M)$ has a canonical (up to a contractible space of choices) quasilocal relative $d$-orientation.
In fact, using the theory of curved Koszul duality from \cite{MR2993002}, one can include the unit and counit and drop oconess.  We will not give this generalization here, since Theorems~\ref{AKSZ construction with boundary} and~\ref{quasilocal frobenius structure}  together suffice to give mapping spaces of the form $\Omega^\bullet(M,Z) \times^h_{\Omega^\bullet(\partial M,Z)} \Omega^\bullet(\partial M,Y)$ $\Pois_0$-structures, where $Z$ is a $\Pois_d$ pointed infinitesimal dg manifold, $Y$ is a coisotropic therein, and $\Omega^\bullet(M,Z) = \Omega^\bullet(M)\otimes Z$ is the complex of $Z$-valued de Rham forms on $M$.  This completes the {Poisson AKSZ construction with coisotropic boundary conditions}.

\subsection*{Proof}

Abstract nonsense of model categories assures that it suffices to check the claim for any particular choice of quasifree resolution.  Let $\LB$ denote the oco dioperad parameterizing Lie bialgebras with all operations in degree $0$.  The dioperad $\Frob_{0,d}^\oco$ is known \cite{MR1960128} to be Koszul with quadratic dual $(\Frob_{0,d}^\oco)^! = \LB\otimes \Frob^\oco_{-1,1-d}$.  We can therefore take $(\mathrm{hFrob}_{0,d}^\oco)^\triangleleft = \bD\bigl( (\LB\otimes \Frob^\oco_{1,1-d})^{\triangleleft,\str}\bigr)$.  After tracking degree shifts, we find that the generators are
$$ G(m,n;m',n') = \begin{cases} 0, & m=n'=0 \\
  \LB(m+n;m'+n')^*[-d(m'+n'-1)+2m+n+m'+2n'-4]
,& \text{else}. \end{cases}$$
The generators are ordered by total arity.

Suppose that $P$ is a quasifree arrowed dioperad and $Q$ is another arrowed dioperad.  One can construct maps $\eta : P \to Q$ by working inductively  in the generators.  The following facts follow from basic obstruction-theoretic yoga (see e.g.\ \cite{me-Qloc}):
\begin{enumerate}
  \item Let $x$ be a closed generator of $P$ of cohomological degree $\deg(x)$ and arity $(m,n;m',n')$.  Suppose that $\H^\bullet Q(m,n;m',n')$ vanishes in cohomological degrees $\bullet < \deg x$.  Then up to a contractible space, the choice of a value of $\eta(x)$ is completely determined by its cohomology class $[\eta(x)] \in \H^{\deg x} Q(m,n;m',n')$.  In particular, if $\H^{\deg x}Q(m,n;m',n') = 0$ as well, then there are \define{no choices}.
  \item Let $x$ be a generator of $P$ of cohomological degree $\deg(x)$ and arity $(m,n;m',n')$ with differential $\partial x$.  Assuming $\eta$ has been defined on generators before $x$, the value of $\eta(\partial x) \in Q(m,n;m',n')$ has already been determined, and is automatically a closed element of cohomological degree $\deg(x) +1$.  Its cohomology class $[\eta(\partial x)] \in \H^{\deg x + 1}Q(m,n;m',n')$ is called the \define{obstruction} for $x$, and $\eta(\partial x) \in Q(m,n;m',n')$ the \define{cochain-level obstruction}.  If earlier choices are changed to homotopy-equivalent choices, then the cochain-level obstruction changes, but its cohomology does not.
  
    Suppose that $\H^\bullet Q(m,n;m',n')$ vanishes in cohomological degrees $\bullet \leq \deg x$.  If the obstruction is non-zero, the construction fails: there is no way to define $\eta(x)$.  If the obstruction vanishes, then the space of choices for $\eta(x)$ is contractible.  In particular, only if $\H^{\deg x + 1}Q(m,n;m',n') \neq 0$ is there a \define{potential obstruction}.
\end{enumerate}

Take $P = \bD\bigl( (\LB\otimes \Frob^\oco_{1,1-d})^{\triangleleft,\str}\bigr)$ and $Q = \Qloc(M^\triangleleft)$.  We have
$$ \H^\bullet\Qloc(M^\triangleleft)(m,n;m',n') \cong \H^\bullet(M;L)[-d(m'+n'-1)]$$
for some submanifold $L$ of $M$.  In particular, it is supported only in cohomological degrees $\bullet \geq d(m'+n'-1)$.  On the other hand, a generator $x$ of arity $(m,n;m',n')$ is in cohomological degree $\deg x = d(m'+n'-1) + (4-2m+n+m'+2n')$.  Since at least one of $m$ or $n'$ is non-zero and since the total arity of any generator is at least three, $2m+n+m'+2n' \geq 4$.
It follows that:
\begin{enumerate}
  \item There are no choices unless $(m,n;m',n')$ is one of the four cases $(0,2;0,1)$, $(1,0;2,0)$, $(1,1;1,0)$, or $(0,1;1,1)$.  These generators are all closed.  Their cohomology classes act on $\End(\H^\bullet(M;\partial M) \to \H^\bullet(M))$ to make $\H^\bullet(M;\partial M)$ into a coalgebra, $\H^\bullet(M)$ into an algebra, $\H^\bullet(M;\partial M)$ into an $\H^\bullet(M)$-module, and $\H^\bullet(M)$ into an $\H^\bullet(M;\partial M)$-comodule.  Since these data are determined by the usual Poincar\'e duality, we see that the cohomology classes $[\eta(x)]$ for the four closed generators are determined up a contractible space of choices.
  \item There are no potential obstructions unless $(m,n;m',n')$ is one of: $(1,1;0,1)$, $(1,0;1,1)$, $(2,0;1,0)$, $(0,1;0,2)$, $(0,3;0,1)$, $(1,0;3,0)$, $(1,2;1,0)$, $(0,1;2,1)$, $(0,2;1,1)$, or $(1,1;2,0)$.  Arguments almost exactly the same as the ones given in~\cite{me-Qloc} imply that the obstructions vanishes. 
  
   Indeed, we may if we so please choose to let the generators of arity $(0,2;0,1)$ and $(1,1;1,0)$ both act as the wedge products.  Then the cochain-level obstructions for the generators of arity $(1,1;0,1)$, $(2,0;1,0)$, $(0,3;0,1)$, and $(1,2;1,0)$ vanish identically.  It follows that the obstructions for these generators vanish for any choices of the closed generators.  On the other hand, as in~\cite{me-Qloc}, we can choose to represent the closed generators by operations with smooth integral kernel, and then the same integral kernels can be used to represent the all closed generators of total arity $3$.  In this case, the cochain-level obstructions for all the generators with potential obstructions are operations with smooth integral kernel, and this kernel only depends on the total arity.  Since the obstruction is known to vanish for the many-to-one generators, it vanishes for all of the generators. \qed
\end{enumerate}

\section*{Acknowledgements}

{This work was supported in part by the NSF grant DMS-1304054. Research at Perimeter Institute for Theoretical Physics is supported by the Government of Canada through the Department of Innovation, Science and Economic Development Canada and by the Province of Ontario through the Ministry of Research, Innovation and Science.}

%\bibliography{ReferencesWithLinks}{}
%\bibliographystyle{alpha}

\end{document}